\documentclass[12pt,reqno,a4paper]{amsart}
\newtheorem{thm}{Theorem}[section]
\newtheorem{defi}{Definition}[section]

\newtheorem{cor}{Corollary}[section]
\newtheorem{pr}{Proposition}[section]
\theoremstyle{definition}
\newtheorem{rem}{Remark}[section]

\newcommand{\be}{\begin{equation}}
\newcommand{\ee}{\end{equation}}
\newcommand{\bea}{\begin{eqnarray}}
\newcommand{\eea}{\end{eqnarray}}
\newcommand{\beb}{\begin{eqnarray*}}
\newcommand{\eeb}{\end{eqnarray*}}
\usepackage{amssymb,amsfonts,amsthm,setspace,indentfirst,mathrsfs}
\usepackage{xcolor}
\usepackage{minibox}
\usepackage{matlab-prettifier}

\usepackage{enumitem}
\usepackage{cite}
\usepackage{float}
\numberwithin{equation}{section}
\begin{document}
%
\title[Pseudosymmetry, Ricci soliton and Curvature Inheritance symmetries of FLRW spacetime]{\bf{Pseudosymmetry, Ricci soliton and Curvature Inheritance symmetries of Friedmann 
Lemaître Robertson Walker spacetime}}

\author[A. A. Shaikh and Kamiruzzaman]{Absos Ali Shaikh$^{* 1}$ and Kamiruzzaman$^{2}$}
\date{\today}

\address{\noindent$^{1,2}$ Department of Mathematics,
	\newline The University of Burdwan, Golapbag,
	\newline Burdwan-713104, West Bengal, India}
\email{aashaikh@math.buruniv.ac.in$^{1}$ ; aask2003@yahoo.co.in$^{1}$ }
\email{kamiruzzaman8145@gmail.com$^{2}$}

\dedicatory{}

\begin{abstract}	
The Friedmann--Lema\^{i}tre--Robertson--Walker (FLRW) spacetime, which was first proposed by Friedmann (1922--1924) and Lema\^{i}tre (1927) and subsequently developed by Robertson and Walker (1935), is an isotropic and homogeneous cosmological model of the universe. This paper addresses a significant gap in the differential geometry literature by providing the comprehensive examination of the curvature properties of the FLRW spacetime. It is demonstrated that the FLRW spacetime satisfies the curvature condition $R \cdot R - Q(S,R)=L_C Q(g,C)$ alongside several pseudosymmetric-type conditions related to the conformal and conharmonic curvature tensors. Furthermore, the Tachibana tensors $Q(g,C)$ and $Q(S,C)$ are found to exhibit a linear dependence on the tensor $(C \cdot R + R \cdot C)$. Additionally, the spacetime is shown to be a 2-quasi-Einstein manifold, generalized Roter type and $\mathrm{Ein}(3)$. The Ricci tensor is shown to be neither cyclic parallel nor of Codazzi type, yet it satisfies several compatibility requirements concerning the Riemann, conformal, projective, conharmonic, and concircular curvature tensors. A thorough analysis of Ricci solitons and curvature inheritance properties reveals that the spacetime admits almost Ricci soliton and $\eta$-Ricci Yamabe soliton structures with respect to the non-Killing soliton vector fields $\frac{\partial}{\partial t}$ and $\frac{\partial}{\partial r}$, but not with respect to $\frac{\partial}{\partial \theta}$. Moreover, the spacetime admits generalized curvature inheritance symmetry properties for the Riemann curvature tensor, as well as for the Weyl conformal, concircular, and conharmonic curvature tensors with respect to the coordinate vector field $\frac{\partial}{\partial t}$ and the gradient of $t$. Later, a comparison of the FLRW and Lema\^{i}tre--Tolman--Bondi (LTB) spacetimes is provided in terms of various curvature-related geometric properties and physical characteristics. Finally, a noteworthy conclusion of the entire study is presented.
	
\end{abstract}

%

\noindent\footnotetext{ $^*$ Corresponding author(Absos Ali Shaikh, E-mail address: aashaikh@math.buruniv.ac.in, aask2003@yahoo.co.in).\\
$\mathbf{2020}$\hspace{5pt}Mathematics\; Subject\; Classification: 53B20; 53B30; 53B50; 53C15; 53C25; 53C35; 83C15. }
\maketitle
%

\section{\bf Introduction}\label{intro}

Let us consider $(M,g_{qs})$ be an $n$-dimensional manifold $(n \ge 3)$, where $M$ is an open and connected subset of $\mathbb{R}^n$, endowed with a semi-Riemannian metric $g_{qs}$ of signature $(\psi,n)$ or $(n,\psi)$, with $0 \le \psi \le n-1$. The metric is Lorentzian when $\psi = 1$ and Riemannian when $\psi = 0$. Throughout this paper, the term spacetime refers to a connected four-dimensional Lorentzian manifold. We denote the Riemann curvature tensor, Ricci tensor, scalar curvature, and Levi-Civita connection by $R_{qust}$, $S_{qs}$, $\kappa$, and $\nabla$, respectively, where the indices $q,u,s,t \in \{1,2,\dots,n\}$. Furthermore, $C_{qust}$, $P_{qust}$, $W_{qust}$, and $K_{qust}$ stand for the conformal, projective, concircular, and conharmonic curvature tensors, respectively. Unless otherwise stated, $M$ $(n \ge 3)$ will always represent a semi-Riemannian manifold.

The energy-momentum tensor (EMT) of type $(0, 2)$ is essential to describe the physical sources in a particular spacetime in the context of general relativity (GR). By affecting quantities like the Ricci tensor, the scalar curvature, and the metric itself, this tensor has a substantial effect on the geometric structure. These components act together to determine the manifold's curvature and general configuration, both of which are crucial for examining the manifold's geometric behavior. Interestingly, the trace of the Ricci tensor, which is derived as the trace of the Riemann curvature tensor, yields the scalar curvature. Therefore, the idea of ``curvature" captures important aspects of spacetime's physical characteristics. For instance, if a Brinkmann wave meets the requirement $R^{ab}_{\ \ qu}R_{abst} = 0$ \cite{Brink1925, SBH21, MS2016, SG1986}, it is considered a pp-wave. Moreover, Petrov type D spacetimes with Weyl tensors that connect pseudosymmetric conditions \cite{DHV2004}. A quasi-Einstein manifold can be regarded as a perfect fluid solution and vice versa, as demonstrated in \cite{C01, SKH11, SYH09}.

In contrast, the structure of a semi-Riemannian manifold is fundamentally determined by its curvature, which depends on the covariant derivatives of tensor fields, possibly of higher order. If the condition $\nabla_a R_{qust}=0$ is satisfied, then the manifold $M$ is called locally symmetric \cite{Cart26}; in this case, the geometry in a neighborhood of each point exhibits symmetry, and the associated local geodesic symmetries act as isometries. This notion constitutes a central concept in differential geometry, reflecting uniform geometric behavior around every point. Subsequently, Ruse \cite{Patt52, Ruse46, Ruse49a, Ruse49b} generalized this framework by introducing the concept of ``kappa'' spaces, which weaken the strict curvature condition of local symmetry and thereby include a broader class of geometric structures. These `kappa' spaces originally arose from investigations of curvature behavior under parallel transport. They allow scalar multiplication of the curvature tensor and are now recognized as a significant class of manifolds. Walker \cite{Walk50} later referred to such manifolds as recurrent manifolds, which stimulated further research and led to the identification of several generalized forms of recurrence. For a comprehensive treatment of these classes of manifolds, the reader may consult \cite{SAR13, SK14, SKA18, SP10, SR11, SRK17, SRK18,SAKI26,ES2024,EDS2022}. Earlier contributions by Shirokov \cite{Shi25} examined manifolds admitting covariantly constant tensor fields, and this line of investigation was subsequently expanded in \cite{Shi98}.

The integrability conditions corresponding to the relation $\nabla R=0$ can be compactly written as $R \cdot R=0$, a formulation that has been thoroughly examined in the classical works of Cartan \cite{Cart46} and Shirokov \cite{Shi25}. A manifold $M$ satisfying $R \cdot R=0$ is called semisymmetric. The more general condition $R \cdot R = L_{R} Q(g,R)$ was systematically investigated by Sinyukov \cite{Sin54}, where $L_{R}$ denotes a smooth function defined on the set $\{x \in M : Q(g,R)x \neq 0\}$ and $Q(g,R)$ is the Tachibana tensor. Manifolds fulfilling this condition are termed as generalized symmetric, and their properties have been discussed in detail \cite{MIKES76, Add1}, with further refined formulations presented in \cite{MIKES88, MIKES92, MIKES96, MSV15, MVH09}.

Over the last several decades, the analysis of symmetry properties in differential geometry and GR has witnessed substantial development, leading to deeper insights into the structure of curved spacetimes. A significant contribution in this direction was provided by Adamów and Deszcz \cite{AD83}, who proposed the notion of second-order symmetry as a natural generalization of classical semisymmetry. Motivated by a careful study of Einstein’s field equations (EFE) and the geometry of totally umbilical hypersurfaces, their results established an important foundation for subsequent investigations of broader symmetric-type conditions.

The classification of spacetimes in GR fundamentally relies on the existence of geometric symmetries, traditionally expressed through the vanishing of the Lie derivative of the metric tensor, i.e., $\pounds_{\zeta} g_{qs} = 0$. While isometries and Killing vector fields (KVFs) provide a robust framework for static and stationary solutions, they are often too restrictive for the description of dynamical, inhomogeneous, or evolving cosmological models. These symmetries are often too rigid for evolving dynamical systems. Consequently, more general symmetry transformations, such as \textit{Ricci inheritance} have gained prominence in the literature \cite{Duggal1999}. This limitation has motivated the study of more general symmetry transformations, most notably, curvature inheritance (CI). A spacetime is said to admit a CI symmetry if the Lie derivative of the Riemann curvature tensor $R^q{}_{ust}$ along a vector field ${\zeta}$ satisfies $\pounds_{\zeta} R^{q}_{ust} = 2\alpha R^q_{ust}$ for some scalar function $\alpha$ \cite{Duggal1992}.

Physically, generalized curvature inheritance (GCI) represents a state where the tidal forces and the gravitational structure evolve proportionally to the background geometry. Unlike curvature collineations (where $\alpha = 0$), which often impose unphysical constraints on the fluid parameters, CI allows for a more flexible coupling between geometry and matter. Through the EFE, this symmetry is inherited by the EMT $T_{qs}$, forcing the matter fields-such as energy density, pressure, and heat flux-to evolve in a self-similar manner \cite{Coley1990}. Such configurations are of particular interest in the study of Ricci solitons (RS), where the metric represents a stationary point of the Ricci flow, effectively modeling universes that expand or contract without altering their underlying topological signatures.

Furthermore, the existence of CI symmetries provides a bridge between the conformal structure of spacetime and its dynamical content\cite{Hall2004}. It ensures that the light-cone structure and the gravitational curvature vary in a synchronized fashion, offering insights into the stability of gravitational vacuums and the potential existence of dark energy-like equations of state ($pressure(p) = -density(\rho)$) in conformally flat spacetimes. Thus, here we mention the physical implications of these generalized symmetries, specifically focusing on their role in defining conservation laws and the thermodynamic evolution of relativistic fluids in non-vacuum backgrounds.

The exploration of spacetime symmetries in GR is fundamental to the classification of exact solutions to the EFE. While the standard approach utilizes isometries characterized by Killing vector fields (KVFs), defined by the vanishing of the Lie derivative of the metric tensor ($\pounds_{\zeta} g_{qs} = 0$), 

A significant geometric feature in this context is the linear dependency between the Lie derivative of the metric tensor and the Ricci tensor. Formally, this is expressed as:
\begin{equation}
    \pounds_{\zeta} S_{qs} = \alpha \pounds_{\zeta} g_{qs} + \beta g_{qs},
\end{equation}
where $S_{qs}$ denotes the Ricci tensor and $\alpha, \beta$ are scalar functions on $M$. The physical significance of this dependency is twofold. Firstly, it represents a state of ``curvature inheritance,'' where the symmetries of the gravitational potential (the metric) are directly mapped onto the matter-energy distribution. As the Ricci tensor is coupled to the EMT $T_{qs}$ via the EFE,
  any linear relation between $\pounds_{\zeta} g$ and $\pounds_{\zeta} R$ imposes stringent constraints on the thermodynamics and equation of state (EoS) of the source fluid \cite{Coley1990}.
Secondly, this linear dependency is the defining characteristic of \textit{Ricci solitons} (RS), which satisfy the equation:
\begin{equation}
    \pounds_{\zeta} g_{qs} + 2S_{qs} + 2\lambda g_{qs} = 0,
\end{equation}
where $\lambda$ is a scalar.
Physically, such manifolds describe spacetimes that evolve purely by self-similarity. They serve as stationary points of the Ricci flow, representing gravitational configurations that expand, shrink, or remain steady without altering their fundamental geometric topology \cite{Hamilton1988}. These linear dependencies play a crucial role in determining the propagation of gravitational degrees of freedom and the resulting conservation laws in non-vacuum spacetimes.

The paper is outlined as follows. In Section~2, we present the historical background of the FLRW spacetime. Section~3 proposes geometrical and physical significances of symmetry and pseudosymmetry properties. Section~4 studies the fundamental definitions and geometric preliminaries required for the study of this spacetime. In Section~5, a detailed investigation of the spacetime is carried out, where several significant results are obtained. Section~6 is devoted to the study of Ricci and Ricci--Yamabe solitons in the context of FLRW spacetime. In Section 7, we examine generalized curvature inheritance and related symmetry properties of the spacetime. Section~8 provides a comparative analysis of the geometric structures associated with the FLRW and LTB spacetimes. Finally, a significant conclusion is presented in section 9.
\section{\bf Historical Background of FLRW spacetime}
The Friedmann-Lema\^{i}tre-Robertson-Walker (FLRW) metric stands as one of the most fundamental achievements in theoretical cosmology. It describes a homogeneous, isotropic expanding or contracting universe and forms the geometric backbone of the standard lambda cold dark matter (briefly, $\Lambda$CDM model) \cite{planck2020}. Its development, spanning from 1917 to the mid-1930s, represents a remarkable confluence of theoretical insight, mathematical rigor, and eventual observational confirmation. 

In 1917, \textbf{Albert Einstein} published \emph{``Kosmologische Betrachtungen zur allgemeinen Relativitätstheorie"} (``Cosmological considerations on the general theory of relativity") \cite{einstein1917}. In this paper, he applied his newly developed general theory of relativity to the universe as a whole. To maintain a static universe-the prevailing belief at the time, he introduced the cosmological constant $\Lambda$ into the field equations:
\begin{equation}
G_{qs} + \Lambda g_{qs} = \frac{8\pi G}{c^4} T_{qs}.
\end{equation}
Although, Einstein later reputedly called this as ``greatest blunder'', the cosmological constant would prove crucial for modern cosmology, particularly in explaining the observed accelerated expansion of the universe \cite{planck2020}.

The first person to derive an expanding universe solution from Einstein's equations was the Russian mathematician and physicist \textbf{Alexander Friedmann}. In 1922, he published \emph{``Über die Krümmung des Raumes"} (``On the curvature of space") in \emph{Zeitschrift für Physik} \cite{friedmann1922}. Friedmann demonstrated that Einstein's field equations admitted non-static solutions, meaning the universe could expand or contract. He derived what are now known as the Friedmann equations:
\begin{align}
\left( \frac{\dot{a}}{a} \right)^2 &= \frac{8\pi G\rho}{3} - \frac{kc^2}{a^2} + \frac{\Lambda c^2}{3}, \
&\frac{\ddot{a}}{a} = -\frac{4\pi G}{3}\left( \rho + \frac{3p}{c^2} \right) + \frac{\Lambda c^2}{3},
\end{align}
where the notations have their usual meaning.
In 1924, Friedmann published a second paper, \emph{``Über die Möglichkeit einer Welt mit konstanter negativer Krümmung des Raumes"} (``On the possibility of a world with constant negative curvature of space"), extending his previous work to include models with negative spatial curvature \cite{friedmann1924}. Einstein initially submitted a brief refutation of Friedmann's 1922 paper to the same journal, claiming to have found a mathematical error. After correspondence with Friedmann, Einstein retracted his objection, acknowledging Friedmann's correct calculations \cite{belenkiy2002}. Tragically, Friedmann died of typhoid in 1925 at the age of 37, never witnessing the observational confirmation of his theoretical predictions.

Working independently and unaware of Friedmann's publications, the Belgian physicist and catholic priest \textbf{Georges Lemaître} arrived at similar conclusions in 1927. His paper \emph{``Un Univers homogène de masse constante et de rayon croissant rendant compte de la vitesse radiale des nébuleuses extra-galactiques"} (``A homogeneous universe of constant mass and increasing radius accounting for the radial velocity of extra-galactic nebulae") was published in the relatively obscure \emph{Annales de la Société Scientifique de Bruxelles} \cite{lemaitre1927}. Lemaître's crucial contribution was not merely re-deriving the expanding universe solutions but linking them to the first observations of galactic redshifts. He estimated a relationship between distance and recession velocity, effectively deriving what would later known as the Hubble-Lemaître law. Lemaître also proposed the concept of a primeval atom or ``Cosmic Egg''—the earliest version of the Big Bang theory \cite{lemaitre1931}. His work remained largely unnoticed until 1930, when \textbf{Arthur Eddington} helped arrange for an english translation, published in the \emph{Monthly Notices of the Royal Astronomical Society} in 1931 \cite{lemaitre1931}. This translation was slightly modified by Lemaître himself, a detail that has led to some historical discussion about the omission of his original calculations \cite{livio2011}.

The theoretical work of Friedmann and Lemaître found its empirical foundation in \textbf{Edwin Hubble}'s 1929 observations. In his paper \emph{``A relation between distance and radial velocity among extra-galactic nebulae"}, Hubble demonstrated a linear relationship between the distance of galaxies and their redshift \cite{hubble1929}. This provided the first strong evidence that the universe was indeed expanding, transforming the FLRW models from mathematical curiosities into viable physical descriptions of our universe. This relationship is now formally known as the Hubble-Lemaître law, a nomenclature officially adopted by the international astronomical union in 2018 to acknowledge Lemaître's earlier theoretical contribution.

The final piece of the puzzle was provided independently by the American mathematician \textbf{Howard P. Robertson} and the British mathematician \textbf{Arthur Geoffrey Walker}. In two papers published in 1935 and 1936- Robertson's \emph{``Kinematics and World-Structure"} \cite{robertson1935, robertson1936} and Walker's \emph{``On Milne's Theory of World-Structure"} \cite{walker1937}—they rigorously proved that the FLRW metric is the \textbf{only} possible metric for a spacetime that is both spatially homogeneous and isotropic. This was a fundamental mathematical proof establishing the geometric foundation for all such models, independent of the specific dynamical equations used to evolve them. Robertson and Walker's work formalized the line element in its modern form as
\begin{equation}
ds^{2} = c^{2}dt^{2} - a^{2}(t)\left[ \frac{dr^{2}}{1-kr^{2}} + r^{2}\left( d\theta^{2} + \sin^{2}\theta d\phi^{2} \right) \right].
\end{equation}
Here, $c$ symbolizes the speed of light in vacuum, $t$ denotes the cosmic time, $r$ represents the radial coordinate, and $a(t)$ is the scale factor describing the time evolution of the universe, indicating whether it expands or contracts with time, such that $a(t)$ is a non-vanishing smooth function of $t$. The geometric structure of the universe is determined by the parameter $k$, where $k=+1$ corresponds to a closed universe with spherical geometry, $k=0$ represents a spatially flat and infinite universe, and $k=-1$ indicates an open universe with hyperbolic geometry.\\
The full four-name acronym \textbf{FLRW} is now widely used to acknowledge all four key contributors \cite{belenkiy2002}. However, different naming conventions persist in various geographical and scientific contexts:

\begin{itemize}
    \item \textbf{FRW} (Friedmann-Robertson-Walker): Common in the United States and parts of Western Europe, often emphasizing the dynamical equations alongside the geometric framework.
    \item \textbf{RW} (Robertson-Walker): Used when focusing specifically on the geometric properties of homogeneity and isotropy, independent of the dynamical evolution \cite{robertson1935, walker1937}.
    \item \textbf{FL} (Friedmann-Lemaître): Often employed to refer specifically to the dynamical solutions—the Friedmann equations and their implications for cosmic evolution \cite{friedmann1922, lemaitre1927}.
\end{itemize}

The centennial of Friedmann's paper in 2022 prompted renewed historical interest with major reviews highlighting the metric's profound impact on cosmology and the rich interconnected history of its discovery \cite{belenkiy2002, livio2011}.

It appears that certain Friedmann--Lema\^{i}tre--Robertson--Walker (FLRW) spacetimes constitute some of the earliest known examples of non-semisymmetric pseudosymmetric warped product manifolds\cite{HV09,{DGJZ-2016}}.
During the mid–twentieth century, the FLRW metric became the central framework of physical cosmology. Early theoretical progress was made by George Gamow and his collaborators \cite{ABG1948}, who applied expanding universe models to explain the primordial abundances of light elements through cosmological nucleosynthesis. Observational confirmation of the hot Big Bang scenario emerged after the discovery of the cosmic microwave background radiation by Penzias and Wilson \cite{PW} in 1965, which provided strong empirical evidence for the thermal relic radiation predicted by FLRW-based cosmology. A major theoretical advance was later introduced by Alan Guth \cite{Guth} in 1981 through the inflationary paradigm, proposing a brief phase of rapid exponential expansion in the early universe that accounts for the observed homogeneity, isotropy, and near-flat spatial geometry while also generating primordial density fluctuations responsible for large-scale structure formation.

At the end of the twentieth century, cosmology experienced a major transformation when observations of distant Type Ia supernovae by Perlmutter et al. \cite{perlmutter1999} and Riess et al. \cite{riess1998} demonstrated that the universe is undergoing accelerated expansion, implying the existence of a dominant component with negative pressure, commonly interpreted as dark energy. High-precision measurements of cosmological parameters were subsequently achieved through observations of the cosmic microwave background by the Planck collaboration \cite{planck2020}, establishing the $\Lambda$CDM model as the standard description of cosmic evolution. More recently, García-Moreno et al. \cite{garciamoreno2026} constructed cosmological spacetimes with time-dependent spatial curvature allowing sign changes and possible topological transitions while maintaining global hyperbolicity. Gomes \cite{gomes2026} developed quasi-dust cosmological models showing that local gravitational potentials in an non-homogeneous universe can produce an apparent dark-energy evolution through backreaction effects. Foidl and Rindler-Daller \cite{foidl2026} introduced a kinematic perspective on late-time cosmic acceleration by studying cosmic expansion within freely falling comoving frames, proposing that the apparent acceleration could emerge from kinematic effects linked to large-scale cosmic structures. Roy et al.\cite{roy2026} investigated photon propagation and Boltzmann dynamics in K-essence cosmology, demonstrating an emergent FLRW geometry with modified causal structure can influence acoustic oscillations and small-scale anisotropies. In addition, Son et al.\cite{son2025} showed that accounting for progenitor age bias in supernova cosmology shifts the inferred value of the Hubble constant toward consistency with BAO and CMB measurements, potentially alleviating the Hubble tension. Historical analyses by Belenkiy \cite{belenkiy2002} and Livio \cite{livio2011} further discussed the development and impact of Friedmann’s cosmological ideas in the emergence of modern cosmology. Although Arslan et al. \cite{ADEHM14} established necessary and sufficient conditions for warped product manifolds, particularly generalized Robertson–Walker (GRW) spacetimes, to satisfy generalized Einstein metric conditions, and Defever et al. \cite{DDHKS00} investigated warped product manifolds with a 1-dimensional base, namely generalized Robertson–Walker spacetimes, the curvature-related geometric structures of FLRW spacetime have not yet been studied comprehensively in the existing literature. In the present paper, we address this gap by establishing several curvature properties, including pseudosymmetry, Roter-type structure, Codazzi-type conditions, and Ricci solitons.

\indent Penrose \cite{Penrose} showed that, in a four-dimensional spacetime, the relation
\[
Q_{[b}C_{a]rs[q}Q_{n]}Q^r Q^s = 0
\]
admits four null solutions, referred to as the principal null directions. If two or more of these null directions coincide, then the Weyl tensor is said to be algebraically special, and the corresponding degeneracy condition is expressed as
\[
Q_{[b}C_{a]rsq}Q^r Q^s = 0.
\]
Such degeneracies give rise to a classification of spacetime into the Petrov types, determined by the multiplicities of the eigenvalues of the self-dual part of the Weyl tensor \cite{MM13}. Moreover, on a Lorentzian manifold, the existence of a null vector $Q$ that is Weyl compatible (or Riemann compatible) implies that the Weyl tensor must be algebraically special. The concept of Weyl compatibility for symmetric tensor fields constitutes a more general requirement than Riemann compatibility. Conditions restricting the form of the Weyl tensor affect both the Petrov classification of the manifold and the properties of its electric and magnetic parts. Various results establish a connection between Weyl compatibility and the Penrose–Debever scheme, showing that a vector field is Weyl compatible precisely when the Weyl tensor is purely electric.
The universal form of the Riemann compatible and Weyl compatible tensor in the FLRW spacetime has been disscused in this paper.



\section{\bf Geometrical and physical significances of symmetry and pseudosymmetry properties}

Symmetry has a fundamental role in the investigation of differential geometry on semi-Riemannian manifolds. Extending the concept of manifolds with constant curvature, the idea of local symmetry was first introduced by Shirokov \cite{Shi25}, and later Cartan \cite{Cart26} provided a complete classification in the Riemannian manifold.

A manifold \(M\) is called locally symmetric when its local geodesic symmetries act as isometries, whereas it is termed globally symmetric if these geodesic symmetries can be extended throughout the entire manifold. Every locally symmetric manifold is globally symmetric; however, the converse does not always hold. For example, if a compact Riemann surface genus greater than one endowed with a usual metric of constant negative curvature \((-1)\), then it is locally symmetric but fails to be globally symmetric \cite{SK14}.

Furthermore, the well-known Cartan--Ambrose--Hicks theorem states that a manifold \(M\) is locally symmetric precisely when $\nabla R = 0.$
In addition, every simply connected and complete locally symmetric manifold is necessarily globally symmetric. Pseudosymmetry generalizes this rigid condition by allowing the curvature tensor to satisfy relaxed, higher-order constraints, which serves as a powerful framework for classifying manifolds with non-uniform curvature. This mathematical relaxation is crucial for identifying structural properties in conformally flat spacetimes, characterizing specific curvature collineations, and modeling spaces like Robertson-Walker geometries. Furthermore, in crystallography, pseudosymmetry allows for the analysis of structural distortions where a low-symmetry crystal lattice mimics a higher-symmetry network, a phenomenon essential for indexing electron backscatter diffraction data and identifying structural pseudosymmetry in atomic configurations.\\
\indent Deszcz tensors are a distinguished class of tensors on the manifold that arise from a thorough analysis of the pseudosymmetry requirement. Eisenhart's work established the Tachibana tensor $Q(g, R)$, which is essential for characterizing pseudosymmetry collineations (PSCs). For the motion group of an $n$-dimensional manifold $M$ to reach its maximal dimension $n(n+1)/2$, it was found that the condition $Q(g, R) = 0$ is both necessary and sufficient. This property is typical of manifolds with constant sectional curvature. Furthermore, the quantity $Q(g,R)_{abcdij}u^{a}v^{b}v^{c}u^{d}x^{i}y^{j}$ admits a geometric interpretation: it is the first-order variation in the sectional curvature of the plane created by the vectors $u$ and $v$ as a result of an infinitesimal rotation of their projections into the plane spanned by $x$ and $y$.\\
\indent The scalar function $L_R$ in a pseudosymmetric manifold coincides with the double sectional curvature and can be geometrically interpreted as the difference between the lengths of two geodesics that correspond to two Levi-Civita parallelogramoids connected by parallel transport around an infinitesimal coordinate parallelogram. Furthermore, a pseudosymmetric manifold is created by projectively transforming a Riemannian manifold onto a semisymmetric manifold. Similarly, a manifold must be pseudosymmetric if it can be projectively mapped onto a pseudosymmetric space.\\
\indent Ricci curvature represents curvature in a specific direction, whereas sectional curvature describes a plane's geometry. Vacuum solutions correspond perfectly to spacetimes where the Ricci curvature vanishes in all directions. Models with a cosmological constant, in which the Ricci curvature is constant in all directions, are the most basic non-vacuum models. This idea also applies to manifolds where parallel transport along a curve or around an infinitesimal coordinate parallelogram does not change the Ricci curvature of each direction. Ricci pseudosymmetric spacetimes are the most basic manifolds for which this invariance is broken. While all pseudosymmetric manifolds are inherently Ricci pseudosymmetric, the opposite is not always true.\\
\indent Geometrically speaking, semisymmetric spaces are manifolds whose sectional curvature does not change up to second order when a plane is parallel transported around a parallelogram of infinitesimal coordinates. On the other hand, pseudosymmetric spaces are distinguished by the fact that their double sectional curvature is isotropic. This means that while a scalar curvature function that depends on a point and two tangent planes at that point can be defined on any manifold, for pseudosymmetric manifolds, this function only varies with the point itself. Subject to further restrictions on the Ricci tensor, an algebraic classification of pseudosymmetric spacetimes reveals that they must either be conformally flat or belong to Petrov type D or N.

In physics, exact symmetries dictate fundamental conservation laws via Noether's theorem, such as space-translation invariance yielding momentum conservation, while also enabling the dimensional reduction of mechanical equations of motion through frameworks like Euler-Poincaré reduction. Pseudosymmetry and approximate symmetries find their primary utility in general relativity and cosmology, where they are used to solve Einstein's field equations for realistic, non-idealized spacetimes, including Vaidya-Bonner-de Sitter and generalized Robertson-Walker cosmological models. Additionally, the subtle breaking of exact symmetry into pseudosymmetry heavily influences condensed matter physics, directly governing macroscopic material behaviors such as piezoelectricity, optical activity, and second-harmonic generation in crystals where minute structural deviations drastically alter electromagnetic properties.

Particularly in molecular systems and crystals, pseudosymmetry serves as an indicator and a precursor to structural phase transitions, which is crucial for the understanding of instabilities \cite{TU25}. In a high symmetry configuration, the existence of pseudosymmetry frequently indicates that the system is not in its most stable form and is vulnerable to a distortion that reduces its symmetry to a more stable configuration. For instance, pseudosymmetry plays a key role in the main theoretical framework known as the Pseudo-Jahn-Teller Effect (PJTE). In high symmetry configurations of polyatomic systems in non-degenerate electronic states, the PJTE characterizes spontaneous distortion and symmetry breaking \cite{BBPE2002}.\\
\indent The ground electronic state's vibronic coupling (mixing) with one or more low-lying excited states with distinct symmetries causes the instability. 
In the end, the connection results in a stable but less symmetrical structure by reducing the curvature of the energy surface in the direction of the distortion.
 In high symmetry molecular systems, the PJTE is thought to be the sole cause of instability and spontaneous symmetry breaking.\\
 \indent Pseudosymmetry is utilized in material science to find possible molecule ferroelectrics. At a particular temperature (the Curie temperature), ferroelectric materials go through a phase transition from a high-symmetry paraelectric phase to a lower-symmetry polar phase. Pseudosymmetry analysis aids in the identification of materials having many reorientable domains and the potential for such transitions\cite{LVGJ}.\\
 \indent Proteins in biological systems display pseudosymmetry due to the imperfect similarity of their subunits. Because more flexible regions may be correlated with larger departures from perfect symmetry, analyzing the pseudosymmetry can shed light on protein mobility, flexibility, and assembly dynamics\cite{Robert15}.\\
 \indent Pseudosymmetry can result from metastable conformational states or particular structural configurations (such as the three-fold pseudosymmetry of the PSI3-IsiA18 supercomplex) in biological systems such as photosynthetic complexes. These structural elements affect the emission and transport of energy.
 Pseudasymmetry can make it more difficult to determine and refine a structure in crystallography.\\
 \indent In a system containing several emitters (such atoms or molecules), pseudosymmetry results in a deviation from the perfect, total symmetry necessary for maximum superradiance. Its main functions include generating subradiant states, altering the collective emission rate, and facilitating phenomena like dynamical symmetry breaking\cite{SHC2022}.\\
 \indent For instance, the traditional Dicke model of superradiance postulates that all emitters interact identically in an area much smaller than the wavelength of light, resulting in a single, extremely superradiant (totally symmetric) state. The term "pseudosymmetry" describes faulty structural or interactional features that resemble this symmetry (such as in extended atomic clouds, disordered arrays, or particular molecular aggregates)\cite{SHC2022}.\\
 \indent Only the symmetric states are occupied during decay in a totally symmetric system. This stringent selection constraint is broken by pseudosymmetry, making it possible to access non-symmetric states (subradiant states). These subradiant states can be long-lived "trapped" states that can be utilized for quantum information storage, and their decay rates are slower than the spontaneous emission rate of a single atom.\\ 
 \indent Directionally dependent (anisotropic) light emission may result from pseudosymmetry's disruption of perfect symmetry. For example, charging external forces can break the mirror symmetry of superradiance in a Bose-Einstein condensate with dipolar interactions, giving control over the asymmetry of superradiant peaks\cite{SHC2024}.\\
 \indent Pseudosymmetry can promote spontaneous as dynamical mirror symmetry breaking in specific situations, such as atoms coupled to a one-dimensional wave path. Large, shot-to-shot variations in the quantity of photons released in various directions are the outcome, and this is indicative of correlated photon state generation\cite{XMS}.\\
 \indent In the circumstances of quasi-normal modes (QNMs), pseudosymmetry is essential to comprehending and describing the spectral instability of the non-Hermitian operators governing black hole disturbances. This characteristic affects the detectability and physical importance of some QNM overtones, especially within the context of the holographic duality{\cite{JMS}}.\\
 \indent The eigenvalues of non-self-adjoint or non-Hermitian operators are called QNM. The system lacks perfect symmetry due to the absence of full self-adjointness, which causes spectral instability{\cite{JMS}}. This implies that the eigen values (QNM frequencies) can drastically shift in response to a minor perturbation of the operator, which represents a slight modification of the black hole's environment or the physical universe.\\
 \indent While the fundamental (longest-lived) QNM mode may be stable and dominate the late-time ringdown signal of a black hole merger, more damped overtones may be extremely sensitive to perturbations, according to the spectral instability identified by pseudosymmetry analysis. This sensitivity leads to complicated transient dynamics and may suggest that a full-time evolution of a perturbation cannot be adequately described by a straightforward sum over QNM.\\
 \indent Pseudosymmetry and the ensuing spectral instability are important for evaluating how well holographic models capture the spectra of real-world tightly coupled systems such as quark-gluon plasma in the setting of gauge/gravity duality. It implies that particular spectrum predictions may be susceptible to tiny modeling error, even when some generic features (such as the ratio of shear viscosity to entropy density) may be stable.\\
 \indent Pseudosymmetry essentially serves as a warning that the seemingly discrete spectrum of QNM may not be a stable, trustworthy set of physical observations under realistic small perturbations. As a result, researchers are compelled to employ more robust tools, such as the pseudospectrum, in order to comprehend the true behavior of these systems \cite{AGL}.
 \indent The mathematical representation of gravitational energy and the phenomenological modes of dark matter and early universe physics are two different and important applications of the concept of pseudosymmetry in gravitational wave physics.\\
 \indent Pseudo-Nambu-Goldstone Boson (PNGB) dark matter arises when a spontaneously broken symmetry (such as a U(1) gauge symmetry) is also intentionally broken in some particle physics models beyond the standard model. Cosmic strings are topological defects in spacetime that can be produced by this high-scale symmetry breaking. As these strings fluctuate and form loops, they create a network that emits a stochastic background of gravitational waves. Using present and upcoming gravitational wave detectors such as LISA, the Big Bang Observer, and pulsar timing arrays, these dark matter theories can be constrained by the characteristics of this gravitational wave background, which are directly related to the energy scale of the pseudosymmetry breaking.\\
 \indent A strong first order phase transition can be linked to the spontaneous symmetry breaking in the early cosmos that results in PNGBs. Future observatories like LISA may be able to detect a stochastic background of gravitational waves produced by this phase shift, offering a novel and complementary method of testing these dark matter models\cite{WNZZ}. As a result, pseudosymmetry plays a variety of roles in gravitational waves, from a crucial theoretical element in contemporary physical models that forecast particular gravitational wave signals from the early universe to a basic mathematical tool for describing energy in general relativity \cite{DDVV}.\\
\section{\bf Preliminaries: Pseudosymmetry, Ricci soliton and curvature inheritance symmetries of semi-Riemannian manifolds}
This section presents a review of the fundamental curvature concepts required for evaluating the various curvature tensors associated with the FLRW metric $(1.1)$. In addition, we outline several significant geometric structures and curvature conditions that possess clear geometric interpretations and are useful in characterizing the symmetry properties of the FLRW spacetime. To identify which of these structures and conditions are satisfied or not satisfied by the metric, it is necessary to examine the relevant curvature-based geometric frameworks.\\
Consider two symmetric tensors of the $(0,2)$-type, $\mathcal{D}$ and $\mathcal{H}$. Their Kulkarni-Nomizu product $(\mathcal{D} \wedge \mathcal{H})$ is defined as
(\cite{DGHS11, DHJKS14, Glog02, G08, Kowa06}):

$$(\mathcal{D}\wedge \mathcal{H})_{qust}=2\mathcal{D}_{q[t}\mathcal{H}_{s]u} +2\mathcal{D}_{u[s}\mathcal{H}_{t]q},$$
here, the symbol $[.]$ denotes antisymmetrization over the specified indices, ensuring that the resulting expression is antisymmetric in those index pairs.

The Christoffel symbols of the second kind $(\Gamma^q_{us})$ with a coordinate basis $(x^1, x^2, \ldots, x^n)$ and a metric tensor $g_{qs}$ are defined as follows:\\
$$
\Gamma^s_{qu} = \frac{1}{2} g^{st} \left( \frac{\partial g_{ut}}{\partial x^q} + \frac{\partial g_{qt}}{\partial x^u} - \frac{\partial g_{qu}}{\partial x^t} \right).
$$

The following are the definitions of $(1,3)$-type curvature tensors, which include the Riemann, conharmonic, concircular, conformal, and projective curvatures, respectively:
\bea
R^h_{qus}&=& 2\left(\Gamma^v_{q[u}\Gamma^h_{s]v} + \partial_{[s}\Gamma^h_{u]q}\right),\nonumber\\
K^h_{qus}&=& R^h_{qus} - \frac{2}{n-2}\left( S^h_{[q}g_{u]s} + \delta^h_{[q}S_{u]s}\right),\nonumber \\
W^h_{qus}&=& R^h_{qus} - \kappa\frac{2}{n(n-1)}\delta^h_{[q}g_{u]s}, \nonumber \\
C^h_{qus}&=& R^h_{qus}+\frac{2}{n-2}\left( \delta^h_{[q}S_{u]s} +S^h_{[q}g_{u]s}\right) -\kappa\frac{2}{(n-1)(n-2)}\delta^h_{[q}g_{u]s} ,\nonumber \\
P^h_{qus}&=& R^h_{qus} -\frac{2}{n-1}\delta^h_{[q}S_{u]s}, \nonumber
\eea
here, $\Gamma^q_{us}$ denotes the connection coefficients, $S^q_u$ implies the Ricci tensor of type $(1,1)$, and $\partial_h$ represents the partial derivative operator $\frac{\partial}{\partial x^h}$.

By reducing indices with the help of the metric tensor $g_{hq}$, the tensors of rank $(0,4)$, namely $R_{hqus},$ $K_{hqus},$ $W_{hqus},$ $C_{hqus},$ and $P_{hqus}$, are obtained. Specifically, they are provided by:

\bea
R_{hqus}&=& g_{h\alpha}(\partial_s \Gamma^\alpha_{qu}-\partial_u \Gamma^\alpha_{qs}+\Gamma^\beta_{qu}\Gamma^\alpha_{\beta s}-\Gamma^\beta_{qs}\Gamma^\alpha_{\beta u}) , \nonumber \\
K_{hqus}&=& R_{hqus} - \frac{1}{n-2} (g\wedge S)_{hqus} , \nonumber \\
W_{hqus}&=& R_{hqus} - \frac{\kappa}{2n(n-1)} (g\wedge g)_{hqus}, \nonumber \\
C_{hqus}&=& R_{hqus}-\frac{1}{n-2}(g\wedge S)_{hqus}+\frac{\kappa}{2(n-1)(n-2)}(g\wedge g)_{hqus} ,\nonumber \\
P_{hqus}&=& R_{hqus} -\frac{1}{n-1}(g_{hs}S_{qu}-g_{qs}S_{hu}). \nonumber
\eea

Let $A$ be a tensor of type $(0,e)$ with $e \geq 1$. Following the notations $(U \cdot A)$ (see \cite{DG02, DGHS98, DH03, SDHJK15, SK14}) and $Q(\mathscr{Y}, A)$ (see \cite{DGPSS11, SDHJK15, SK14, Tach74}), the associated tensor of type $(0,e+2)$ are defined as follows:
\begin{align*}
(U\cdot A)_{b_1b_2\cdots b_e ws} &= -\left[ U^\alpha_{wsb_1}A_{\alpha b_2\cdots b_e}+ \cdots + U^\alpha_{wsb_e}A_{b_1\cdots b_{e-1}\alpha}\right], \\
Q(\mathscr{Y},A)_{b_1b_2\cdots b_e qu} &= \mathscr{Y}_{ub_1}A_{qb_2\cdots b_e}+ \cdots + \mathscr{Y}_{ub_e}A_{b_1b_2\cdots q} \\
&\quad - \mathscr{Y}_{qb_1}A_{ub_2\cdots b_e}- \cdots - \mathscr{Y}_{qb_e}A_{b_1b_2\cdots u},
\end{align*}
where $\mathscr{Y}_{qu}$ is a symmetric tensor of type $(0,2)$ and $U^h_{qus}$ is a tensor of type $(1,3)$.
\begin{defi} $($\cite{AD83, Cart46, Chak87, Chak88, Desz92, Desz93, DGHZ15, DGHZ16, SK14, SKppsnw, Szab82, Szab84, Szab85}$)$ 

A manifold $M$ is said to be $A$-pseudosymmetric with respect to $U$ if the product $U \cdot A$ is proportional to $Q(g,A)$; that is,
\[
(U\cdot A)_{b_1b_2\cdots b_e ws} = f_A \, Q(g,A)_{b_1b_2\cdots b_e ws},
\]
where $f_A$ is a smooth function on $M$.

Furthermore, if
\[
(U \cdot A)_{b_1b_2\cdots b_e ws} = \bar{f}_A \, Q(S,A)_{b_1b_2\cdots b_e ws},
\]
for some smooth function $\bar{f}_A$ defined on $M$, then $M$ is called Ricci-generalized $A$-pseudosymmetric with respect to $U$.
\end{defi} 

The manifold is denoted as pseudosymmetric or Ricci-generalized pseudosymmetric (RGPS) whenever both $U$ and $A$ are chosen as the Riemann curvature tensor $R$. It is possible to classify different types of pseudosymmetric and RGPS manifolds by looking at different curvature tensors, such as $R,$ $S,$ $C,$ $W,$ $K$ and $P.$

If $U \cdot A = 0$ is hold, then a manifold is $A$-semisymmetric with respect to a tensor $U$. Specifically, $R \cdot R = 0$ holds for a semisymmetric manifold \cite{Cart46, Szab82, Szab84, Szab85}. Within the more general category of pseudosymmetric manifolds, semisymmetric manifolds constitute a special class that illustrates the importance of particular curvature relations in differential geometry. This categorization demonstrates how specific symmetric and pseudosymmetric qualities in the geometric structure of manifolds are defined by criteria such as $U \cdot A = 0$.

\begin{defi}$($\cite{C01,DGHZ16, DGJZ-2016, DGP-TV-2011, S09, SKH11, SK19,SYH09}$)$ 
 
A manifold $M$ is referred to as an $m$-quasi-Einstein manifold if the tensor $(S - \zeta g)$ has rank $m$ for some scalar $\zeta$ and $0 \leq m \leq (n-1)$. Specifically, when $m = 1$ (or $m = 0$, respectively), the manifold is called as quasi-Einstein \cite{SKH11, SYH09} (or Einstein, respectively). \\
Moreover, if $\zeta = 0$ in the quasi-Einstein condition, the manifold is signified as Ricci simple.
 
\end{defi}

Notably, the Kaigorodov spacetime \cite{SDKC19} is an example of an Einstein manifold, but the Robertson–Walker spacetime \cite{ARS95, Neill83, SKMHH03} is a quasi-Einstein manifold. Examples of 2-quasi-Einstein manifolds are the point-like global monopole spacetime \cite{SASK_pgm_2024} and the Som–Raychaudhuri spacetime \cite{SK16srs}. Additionally, spacetimes like Morris–Thorne \cite{ECS22}, Gödel \cite{DHJKS14}, and Vaidya \cite{SKS19} are categorized as Ricci simple manifolds.

\begin{defi} (\cite{Gray78, SB08, F81, S81})

A manifold $M$ is said to possess a \emph{cyclic parallel Ricci tensor} if at every point of $M$, the Ricci tensor $S_{qs}$ satisfies
$$\nabla_hS_{qs}+\nabla_qS_{hs}+\nabla_sS_{hq}=0.$$
Furthermore, the Ricci tensor $S$ is said to be of \emph{Codazzi type} if it fulfills the condition

\end{defi}
The Ricci tensor in the $(t\text{-}z)$-plane wave spacetime \cite{EC21} is of the Codazzi type, which should be noted, while cyclic parallel Ricci tensor behavior has been noted in the Gödel spacetime \cite{DHJKS14}.

\begin{defi} $($\cite{Bess87, SK14, SK19}$)$ 

A manifold $M$ is called an Einstein manifold of level $4$ provided that it satisfies the relation
\[
\bar{\Pi}_1 S^4_{qs} + \bar{\Pi}_2 S^3_{qs} + \bar{\Pi}_3 S^2_{qs} + \bar{\Pi}_4 S_{qs} + \bar{\Pi}_5 g_{qs} = 0,
\quad (\bar{\Pi}_1 \neq 0),
\]
where the coefficients $\bar{\Pi}_i$ $(i=1,\dots,5)$ are smooth functions on $M$.

If $\bar{\Pi}_1 = 0$ while $\bar{\Pi}_2 \neq 0$, then $M$ is termed an Einstein manifold of level $3$. Likewise, for $\bar{\Pi}_1 = \bar{\Pi}_2 = 0$ and $\bar{\Pi}_3 \neq 0$, the manifold is called an Einstein manifold of level $2$.
\end{defi}

It has been observed that the Vaidya–Bonner spacetime \cite{SDC} and the Lifshitz spacetime \cite{SSC19} serve as examples of level 3 Einstein manifolds, while the Siklos spacetime \cite{SDKC19}, Lemos black hole spacetime \cite{SASK_LBH_2025} and the Reissner-Nordstr\"{o}m-de Sitter spacetime \cite{SK2026} represent level 2 Einstein manifolds.

\begin{defi} (\cite{DGJPZ13, DGJZ-2016, SK19})

A manifold $M$ is said to be of generalized Roter type (GRT) if its Riemann curvature tensor $R$ admits a decomposition of the form
\[
R_{qust} = (\mathcal{Z}_{11} S^2_{qs} + \mathcal{Z}_{12} S_{qs} + \mathcal{Z}_{13} g_{qs}) \wedge S^2_{ut}
    + (\mathcal{Z}_{22} S_{qs} + \mathcal{Z}_{23} g_{qs}) \wedge S_{ut}
    + \mathcal{Z}_{33} (g \wedge g)_{qust},
\]
where each $\mathcal{Z}_{ij}$ is a smooth scalar function on $M$.

Moreover, if the coefficients $\mathcal{Z}_{ij}$ are chosen in such a way that $R$ reduces to a linear combination consisting solely of the terms $(g \wedge g)_{qust}$, $(g \wedge S)_{qust}$, and $(S \wedge S)_{qust}$, then $M$ is called a \emph{Roter type} manifold \cite{DG02, DGP-TV-2011}.

\end{defi}

We note that the Reissner-Nordstr\"{o}m-de Sitter spacetime \cite{SK2026}, Melvin magnetic spacetime \cite{SAAC20}, and Robinson–Trautman spacetime \cite{SAA18} are identified as Roter type manifolds. In contrast, the Lifshitz spacetime \cite{SSC19} and Vaidya–Bonner spacetime \cite{SDC} are examples of GRT manifolds.

\begin{defi}$($\cite{DD91,DGJPZ13, MM12a, MM12b, MM13,MM22b}$)$
 
Let $D$ be a \((0,4)\)-tensor on $M$. The Ricci tensor of $M$ is known as $D$-compatible if the following cyclic identity is satisfied on \(M\):\[\underset{q,u,s}{\mathcal{S}}\, S_q^{\,w}\,D_{ustw}=0,\]where \(\mathcal{S}\) represents the cyclic sum over the indices \(q,u,s\). Moreover, a \(1\)-form \(\eta\) is called $D$-compatible whenever the tensor \(\eta \otimes \eta\) fulfills the above $D$-compatibility condition.
\end{defi}

Similarly we can define the compatibility with other curvature tensor by replacing the tensor $D.$ 

%

\begin{defi} (\cite{Yano, SASK_LBH_2025})

A vector field $\zeta$ on $M$ is said to be a \emph{Killing vector field} if it satisfies
\[
\pounds_{\zeta} g_{qs} = 0,
\]
where $\pounds_{\zeta}$ denotes the Lie derivative taken along the direction of $\zeta$.

\end{defi}
In 1988, Hamilton \cite{Hamilton1988} introduced the Yamabe flow as an analogue of the Ricci flow. Later, G\"{u}ler and Cr\u{a}\c{s}m\u{a}reanu \cite{Guler2019} developed the Ricci--Yamabe flow, which combines the features of both the Ricci and Yamabe flows within a unified geometric evolution framework. The self-similar solutions corresponding to this combined flow are called Ricci--Yamabe solitons, whereas the self-similar solutions associated exclusively with the Yamabe flow are known as Yamabe solitons (YS).

The concept of Ricci--Yamabe solitons broadens the study of geometric flows by establishing a unified framework that incorporates the properties of both Ricci and Yamabe solitons. This approach offers deeper insight into the evolution of manifolds under the combined influence of these two flows. The development of the Ricci--Yamabe flow represents a significant advancement in geometric analysis, extending the pioneering work of Hamilton to new directions of investigation, particularly in the context of semi-Riemannian geometry.

\begin{defi}\label{def_ARS} (\cite{Siddiqi2020})

The manifold $M$ is said to admit an $\eta$-Ricci--Yamabe soliton if there exists a non-zero $1$-form $\eta$, together with constants $\phi_1$, $\phi_2$, and $\mu$, such that the following condition is satisfies:
$$
\frac{1}{2} \pounds_\zeta g_{qs} + \phi_1 S_{qs} + \left(\mu - \frac{1}{2} \phi_2 \kappa \right) g_{qs} + \alpha_1 \eta_q \eta_s = 0,
$$
\end{defi}
\noindent holds, where $\alpha_1$ is a scalar. If the constants $\phi_1$, $\phi_2$, $\mu$, and $\alpha_1$ are extended to smooth functions defined on the manifold, the resulting structure is called an almost $\eta$-Ricci--Yamabe soliton. This generalization provides greater flexibility and facilitates the study of more intricate geometric properties within the broader framework of differential geometry.

In semi-Riemannian manifolds, different selections of the constants $\phi_1$, $\phi_2$, and $\mu$ lead to various types of soliton structures. For instance, when $\phi_1 = 0$ and $\phi_2 = 2$, the resulting structure corresponds to a Yamabe soliton, whereas the choice $\phi_1 = 1$ and $\phi_2 = 0$ yields a Ricci soliton. Such solutions play an essential role in describing the behavior of a manifold under the respective geometric flows. Moreover, if the parameters $\phi_1$, $\phi_2$, and $\mu$ are considered as smooth, non-constant functions defined on the manifold, then the corresponding structure is known as an almost Ricci--Yamabe soliton \cite{Siddiqi2020}. This generalization provides additional flexibility and allows the investigation of more intricate geometric structures.\\
In particular:
\begin{itemize}
\item[(i)] If $\phi_1=1$ and $\phi_2=\alpha_1=0,$ then the manifold is said to admit Ricci soliton\cite{Hamilton1982}.
\item[(ii)] If $\phi_1=1$ and $\phi_2=0,$ then it is reffered to as an $\eta$- Ricci soliton\cite{Cho2009, Pigola2011}.
\item[(iii)] If $\alpha_1=0,$ then it is reffered to as a Ricci--Yamabe soliton\cite{Siddiqi2020, SMM2025}.
\end{itemize}
As a generalization of curvature collineation (CC)\cite{KLD1969, KLD1970} $\pounds_\zeta \widetilde{R} = 0,$ Duggal \cite{Duggal1992} introduced the notion of curvature inheritance (CI) as $
\pounds_\zeta \widetilde{R}^h_{qs} = \mathcal{B}\,\widetilde{R}^h_{qs},
$

\begin{defi}\label{def_GCI} (\cite{ShaikhDatta2022})

A manifold $M$ is said to possess 
generalized curvature inheritance (GCI) with respect to the $(0,4)$-type curvature tensor $R_{qust}$ 
if there exists a non-Killing vector field $\zeta$, such that
\[
\pounds_{\zeta} R_{qust}
= \mathcal{J} R_{qust}
+ \mathcal{J}_1 (g \wedge g)_{qust}
+ \mathcal{J}_2 (g \wedge S)_{qust}
+ \mathcal{J}_3 (S \wedge S)_{qust},
\]
where $\mathcal{J}, \mathcal{J}_1, \mathcal{J}_2,$ and $\mathcal{J}_3$ 
are smooth scalar functions on $M$.
\end{defi}
In particular:
\begin{itemize}
\item[(i)] If $\mathcal{J}_i \neq 0$ and $\mathcal{J}_i = 0$ for $i=1,2,3$, then $M$ 
satisfies the condition of curvature inheritance for $R_{qust}$.

\item[(ii)] If $\mathcal{J} = 0$ and $\mathcal{J}_i = 0$ for $i=1,2,3$, 
then the above relation reduces to the condition of CC 
for the tensor $R_{qust}$.
\end{itemize}

\section{\bf Friedmann–
Lemaître–Robertson–Walker spacetime admitting pseudosymmetries}
%
In the coordinates  $(t,r,\theta,\phi)$, the metric tensor  of FLRW spacetime is given by:
$$g=\left(\begin{array}{cccc}
	c^2 & 0 & 0 & 0 \\
	0 & -\frac{a(t)^2}{1-kr^2} & 0 & 0 \\
	0 & 0 & -a(t)^2 r^2 & 0 \\
	0 & 0 & 0 & -a(t)^2 r^2 \sin^2\theta 
\end{array}
\right).
$$
That is,
 \begin{eqnarray}
g_{11}=c^2,\; g_{22}=-\frac{a(t)^2}{1-kr^2}, \ \ g_{33}=-a(t)^2 r^2, \;g_{44}=-a(t)^2 r^2 \sin^2\theta , g_{qs}=0, \text{ for }  q\neq s \; \label{1}
 \end{eqnarray}
where $c$ is speed of light in vacuum, $t$ is cosmic (proper) time measured by comoving observers, $a(t)$ is scale factor that determine how the universe expand or contracts over time, $k$ is the parameter determine shape of the universe and $r$ is radial coordinates respectively.  


In the sequel, we cosider the expressions written as follows:

$U=-1+k^2,$ $U_1=\dot{a}^2+a\ddot{a},$ $U_2=\dot{a}^2-a\ddot{a},$ $U_3=2\dot{a}^2+a\ddot{a},$ and $U_4=\dot{a}^2+2a\ddot{a},$ 
$$\dot{a}=\frac{d(a(t))}{dt},\text{ } \ddot{a}=\frac{d^2(a(t))}{dt^2},\text{ } a^{(3)}=\frac{d^3(a(t))}{dt^3}.$$

The non-vanishing Christoffel symbols of the second kind ($\Gamma^h_{qs}$), which are stated as follows, are obtained from the computations:

 \begin{eqnarray}\label{2}
 \begin{cases}

 	\Gamma^2_{12}=\frac{\dot{a}}{a}=\Gamma^3_{13}=\Gamma^4_{14}, \ \ 
 	\Gamma^1_{22}=\frac{a\dot{a}}{c^2U}, \\
 	\Gamma^3_{23}=\frac{1}{r}=\Gamma^4_{24}, \ \ 
 	\Gamma^1_{33}=\frac{r^2a\dot{a}}{c^2}=\frac{1}{\sin^2\theta}\Gamma^1_{44},\\
 	\Gamma^2_{33}=rU=\frac{1}{\sin^2\theta}\Gamma^2_{44}, 
 	\Gamma^4_{43}=\cot\theta,\ \ \Gamma^3_{44}=-\cos\theta \sin\theta. \\
\end{cases}
  \end{eqnarray}
 
The calculated non-zero components of the Riemann curvature tensor and the Ricci tensor, which might not vanish are as follows:
\begin{eqnarray}\label{3}
 \begin{cases}

 	R_{1212}=-\frac{a\dot{a}}{U}, \ \ 	
 	R_{1313}=\frac{r^2a\ddot{a}}=\frac{1}{\sin^2\theta} R_{1414}, \\
 	R_{2323}=\frac{r^2a^2\dot{a}^2}{c^2U}=-\frac{1}{\sin^2\theta}R_{2424}, \  \
 	R_{3434}=-\frac{r^4a^2}{c^2}(kc^2+\dot{a}^2) \sin^2\theta; \\
 \end{cases}
  \end{eqnarray}
 
 \begin{eqnarray}\label{4}
 \begin{cases}
 	S_{11}=\frac{3\ddot{a}}{a}, \ \ 
 	S_{22}=\frac{2\dot{a}^2+a\ddot{a}}{c^2U}, \\
 	S_{33}=-\frac{c^2kr^2+r^2(2\dot{a}^2+a\ddot{a})}{c^2}=\frac{1}{\sin^2\theta}S_{44}. \ \ 
  \end{cases}
   \end{eqnarray}

  Finally, the scalar curvature is given by
  \begin{equation}
      \kappa=g_{ab}\,S^{ab}=\frac{2(c^2k+3r^2(\dot{a}^2+a\ddot{a}))}{a^2c^2}.\label{5}
  \end{equation}

Let $\mathcal{F}^{1} = R \cdot R,$ $\mathcal{F}^{2} =\nabla R$ and $\mathcal{M}^1 = Q(S, R).$ Taking into account their inherent symmetries, the non-vanishing components of $\mathcal{F}^{1},$ $\mathcal{F}^{2}$ and $\mathcal{M}^1$ are listed below:

\begin{eqnarray}\label{RR}
\begin{cases}
		\mathcal{F}^1_{1223,12}=\frac{r^2a\ddot{a}U_2}{c^2U}=\frac{1}{\sin^2\theta}\mathcal{F}^1_{1224,12}=-\mathcal{F}^1_{1223,13}=-\frac{1}{\sin^2\theta}\mathcal{F}^1_{1214,14}, \\ 
	\mathcal{F}^1_{1223,13}=\frac{r^4a\ddot{a} \sin^2\theta(c^2k+U_2)}{c^2}=\mathcal{F}^1_{1334,14}, \\ 
	\mathcal{F}^1_{2434,23}=\frac{r^4ka^2\dot{a}^2\sin^2\theta}{c^2U}=-\mathcal{F}^1_{2334,24}; 
	\end{cases}
	\end{eqnarray}

\begin{eqnarray}\label{R}
\begin{cases}
\mathcal{F}^{2}_{1212,1}=\frac{\dot{a}\ddot{a}-aa^{(3)}}{U},\\
\mathcal{F}^{2}_{1223,3}=\frac{r^2a\dot{a}U_2}{c^2U}=\mathcal{F}^{2}_{1323,2}=\frac{1}{\sin^2\theta}\mathcal{F}^{2}_{1424,2}=\frac{1}{\sin^2\theta}\mathcal{F}^{2}_{1224,4}=-\frac{1}{2}\mathcal{F}^{2}_{2323,1}=-\frac{1}{2\sin^2\theta}\mathcal{F}^{2}_{2424,1},\\
\mathcal{F}^{2}_{1313,1}=-r^2(\dot{a}\ddot{a}-aa^{(3)})=\frac{1}{\sin^2\theta}\mathcal{F}^{2}_{1414,1},\\
\mathcal{F}^{2}_{1334,4}=-\frac{r^4a\ddot{a} \sin^2\theta(c^2k+U_2)}{c^2}=-\mathcal{F}^{2}_{1434,3}=\frac{1}{2}\mathcal{F}^{2}_{3434,1},\\
\mathcal{F}^{2}_{2334,4}=-r^3ka^2\sin^2\theta=\mathcal{F}^{2}_{2434,3}=\frac{1}{2}\mathcal{F}^{2}_{3434,2};\\
	\end{cases}
	\end{eqnarray}

\begin{eqnarray}\label{Q(S,R)}
\begin{cases}	
\mathcal{M}^1_{1332,12}=\frac{r^2a\ddot{a}U_2}{c^2U}=\frac{1}{\sin^2\theta}\mathcal{M}^1_{1424,12}=-\mathcal{M}^1_{1223,13}=-\frac{1}{\sin^2\theta}\mathcal{M}^1_{1224,14},\\
	\mathcal{M}^1_{1434,13}=-\frac{r^4a\ddot{a} \sin^2\theta(2c^2k+U_2)}{c^2}                                =-\mathcal{M}^1_{1334,14}, \\
	\mathcal{M}^1_{1213,23}=-\frac{kr^2a\ddot{a}}{U}=\frac{1}{\sin^2\theta}\mathcal{M}^1_{1214,24},\\
	\mathcal{M}^1_{2434,23}=\frac{r^4ka^2\sin^2\theta}{c^2U}=-\mathcal{M}^1_{2334,24}; \\
	\end{cases}
		\end{eqnarray}



%
%
%
%

By considering the symmetry characteristics of the conformal curvature tensor $C_{qust}$, the non-zero components can be expressed as follows.

$$\begin{array}{c}
	C_{1212}=-\frac{c^2k}{3U},\\ 
	C_{1313}=-\frac{c^2kr^2}{6}=\frac{1}{\sin^2\theta}C_{1414}, \\
	C_{2323}=-\frac{kr^2a^2}{6U}=\frac{1}{\sin^2\theta}C_{2424}, \\ 
	C_{3434}=-\frac{kr^4a^2\sin^2\theta}{3}.\\
\end{array}$$
The non-vanishing components of the covariant derivatives of the conformal curvature tensor are given by:
$$\begin{array}{c}
C_{1212,1}=\frac{2c^2k\dot{a}}{3aU},\\
C_{1212,2}=\frac{2c^2k}{3rU},\\
C_{1213,3}=\frac{c^2kr}{2}=\frac{1}{\sin^2\theta}C_{1214,4}=-\frac{2}{3}C_{1313,2}=-\frac{2}{3\sin^2\theta}C_{1414,2},\\
C_{1223,3}=-\frac{kr^2a\dot{a}kr}{2U}=\frac{1}{\sin^2\theta}C_{1224,4}=-\frac{2}{3}C_{2323,1}=-\frac{2}{3\sin^2\theta}C_{2424,1},\\
C_{1313,1}=\frac{c^2kr^2\dot{a}}{3a}=\frac{1}{\sin^2\theta}C_{1414,1},\\
C_{1334,4}=-\frac{kr^4a\dot{a}\sin^2\theta}{2}=-C_{1434,3}=\frac{4}{3}C_{3434,1},\\
C_{2323,2}=-\frac{kra^2}{3U}=-\frac{1}{\sin^2\theta}C_{2424,2},\\
C_{2334,4}=-\frac{kr^3a^2\sin^2\theta}{2}=C_{2434,3}=\frac{4}{3}C_{3434,2}.\\   
\end{array}$$

Let $\mathcal{F}^{4}=R\cdot C$, $\mathcal{F}^{5}=C\cdot R$, $\mathcal{F}^{6}=C\cdot C$ and $\mathcal{M}^{3}=Q(g,C)$, $\mathcal{M}^{4}=Q(S,C)$. 
Taking into account their symmetry properties, the non-vanishing components of $\mathcal{F}^{4}$, $\mathcal{F}^{5}$, $\mathcal{F}^{6}$, $\mathcal{M}^{3}$, and $\mathcal{M}^{4}$ are given below:

\begin{eqnarray}\label{RC}
\begin{cases}
		\mathcal{F}^{4}_{1223,13}=\frac{kr^2a\ddot{a}}{2U}=\frac{1}{\sin^2\theta}\mathcal{F}^{4}_{1224,14}, \\ 
		\mathcal{F}^{4}_{1434,13}=-\frac{kr^4a\ddot{a}\sin^2\theta}{2}=-\mathcal{F}^{4}_{1334,14},\\
		\mathcal{F}^{4}_{1213,23}=-\frac{kr^2\dot{a}^2}{2U}=\frac{1}{\sin^2\theta}\mathcal{F}^{4}_{1214,24},\\
		\mathcal{F}^{4}_{2434,23}=-\frac{kr^4a^2\dot{a}^2\sin^ 2\theta}{2c^2U}=-\mathcal{F}^{4}_{2334,24}; \\
	
	\end{cases}
			\end{eqnarray}

\begin{eqnarray}\label{CR}
\begin{cases}
	\mathcal{F}^{5}_{1323,12}=-\frac{kr^2U_2}{3U}=\frac{1}{\sin^2\theta}\mathcal{F}^{5}_{1424,12}=2\mathcal{F}^{5}_{1223,13}=\frac{2}{\sin^2\theta}\mathcal{F}^{5}_{1224,14}, \\
	\mathcal{F}^{5}_{1434,13}=\frac{kr^4(c^2k+U_2 )\sin^2\theta}{6}=-\mathcal{F}^{5}_{1334,14}, \\
	\mathcal{F}^{5}_{2434,23}=\frac{kr^4a^2\sin^2\theta}{6U}=-\mathcal{F}^{5}_{2334,24}; 
	\end{cases}
			\end{eqnarray}

\begin{eqnarray}\label{CC}
\begin{cases}
\mathcal{F}^6_{1223,13}=\frac{c^2kr^2}{12U}=\frac{1}{\sin^2\theta}\mathcal{F}^{6}_{1224,14}=-\mathcal{F}^{6}_{1213,23}=\frac{1}{\sin^2\theta}\mathcal{F}^{6}_{1214,24}, \\ 
	\mathcal{F}^{6}_{1434,13}=-\frac{c^2kr^4\sin^2\theta}{r^6}=-\mathcal{F}^{6}_{1334,14}, \\
	\mathcal{F}^{6}_{2434,23}=\frac{kr^4a^2\sin^2\theta}{12U}=-\mathcal{F}^{6}_{2334,24};
\end{cases}
		\end{eqnarray}

\begin{eqnarray}\label{Q(g,C)}
\begin{cases}
\mathcal{M}^3_{1223,13}=\frac{c^2kr^2a^2}{2U}=\frac{1}{\sin^2\theta}\mathcal{M}^3_{1224,14}=-\mathcal{M}^3_{1213,23}=\frac{1}{\sin^2\theta}\mathcal{M}^3_{1214,24},\\	
	\mathcal{M}^3_{1434,13}=-\frac{c^2kr^4\sin^2\theta}{r^2}=-\mathcal{M}^3_{1334,14}, \\
	\mathcal{M}^3_{2434,23}=-\frac{kr^4a^4\sin^2\theta}{2U}=-\mathcal{M}^3_{2334,24}; \\
\end{cases}
		\end{eqnarray}

\begin{eqnarray}\label{Q(S,C)}
\begin{cases}
\mathcal{M}^{4}_{1323,12}=-\frac{kr^2U_2}{3U}=\frac{1}{\sin^2\theta}\mathcal{M}^{4}_{1424,12},\\
\mathcal{M}^{4}_{1223,13}=\frac{kr^2(2c^2+(4\dot{a}^2+5a\ddot{a}))}{6U}=\frac{1}{\sin^2\theta}\mathcal{M}^{4}_{1224,14},\\
\mathcal{M}^{4}_{1434,13}=-\frac{kr^4\sin^2\theta(c^2k+(2\dot{a}^2+7a\ddot{a}))}{6}=-\mathcal{M}^{4}_{1334,14}, \\
\mathcal{M}^{4}_{1213,23}=\frac{kr^2(2c^2+3(2\dot{a}^2+a\ddot{a}))}{6U}=\frac{1}{\sin^2\theta}\mathcal{M}^{4}_{1214,24}, \\ 

		\mathcal{M}^{4}_{2434,23}=-\frac{kr^4a^2\sin^2\theta(c^2k+3(2\dot{a}^2+a\ddot{a}))}{6c^2U}=-\mathcal{M}^{4}_{2334,24}. 	
\end{cases}
		\end{eqnarray}

By considering the symmetry properties of the projective curvature tensor $P_{qust}$, its non-vanishing components are presented below:

$$\begin{array}{c}
	P_{1221}=-\frac{2U_2}{3U},  \\
    P_{1331}=\frac{r^2(c^2k+2U_2)}{3}=\frac{1}{\sin^2\theta}P_{1441},\\
	P_{2323}=\frac{r^2a^2U_2}{3c^2U}=\frac{1}{\sin^2\theta}P_{2424},\\
	P_{2332}=-\frac{a^2r^2(c^2k+U_2)}{3c^2U}=\frac{1}{\sin^2\theta}P_{2442}, \\
	P_{3434}=\frac{r^4a^2\sin^2\theta(c^2k+U_2)}{3c^2}=-P_{3443}.
\end{array}$$

By considering the symmetry characteristics of the concircular curvature tensor $W_{qust}$, the non-vanishing components are given below:

$$\begin{array}{c}
W_{1212}=\frac{c^2k+3U}{6U},
 W_{1313}=-\frac{r^2(c^2k+U_2)}{6}=\frac{1}{\sin^2\theta}W_{1414}, \\
	W_{2323}=\frac{a^2r^2(c^2k+3U_2)}{6c^2U}=\frac{1}{\sin^2\theta}W_{2424}, 
	W_{3434}=-\frac{a^2r^4\sin^2\theta(5c^2k+3U_2)}{6c^2}. \\ 
\end{array}$$

By taking into consideration the symmetry properties of the conharmonic curvature tensor $K_{qust}$, the following components are found to be non-vanishing:
	
$$\begin{array}{c}
	K_{1212}=\frac{U_1}{U},\ \ 
	K_{1313}=-\frac{r^2(c^2k+U_1)}{2},\ \
	K_{1414}=\frac{r^2\sin^2\theta(c^2k+2U_1)}{2}, \\
K_{2323}=-\frac{a^2r^2(c^2k+2U_1)}{2c^2U}=\frac{1}{\sin^2\theta}K_{2424}, \ \ 
	K_{3434}=\frac{ar^4\sin^2\theta U_1}{c^2}.
\end{array}$$

Let $\mathcal{F}^{7}=C \cdot K$, $\mathcal{F}^{8}=K \cdot C$, $\mathcal{F}^{9}=K \cdot K$, and $\mathcal{M}^{5}=Q(g,K)$. 
By taking their symmetry properties into account, the non-vanishing components of $\mathcal{F}^{7}$, $\mathcal{F}^{8}$, $\mathcal{F}^{9}$, and $\mathcal{M}^{5}$ are obtained as follows:


\begin{eqnarray}\label{CK}
\begin{cases}
	\mathcal{F}^{7}_{1223,13}=-\frac{c^2kr^2}{12U}=\frac{1}{\sin^2\theta}\mathcal{F}^{7}_{1224,14}=-\mathcal{F}^{7}_{1213,23}=-\frac{1}{\sin^2\theta}\mathcal{F}^{7}_{1214,24}, \\
	\mathcal{F}^{7}_{1434,13}=\frac{c^2kr^4\sin^2\theta}{12}=-\mathcal{F}^{7}_{1334,14}, \\
	\mathcal{F}^{7}_{2434,23}=\frac{kr^4a^2\sin^2\theta}{12U}=-\mathcal{F}^{7}_{2334,24}; 
\end{cases}
		\end{eqnarray}
		
\begin{eqnarray}\label{KC}
\begin{cases}
	\mathcal{F}^{8}_{1223,13}=-\frac{kr^2(c^2k+2U_1)}{4U}=\frac{1}{\sin^2\theta}\mathcal{F}^{8}_{1224,14}=-\mathcal{F}^{8}_{1213,23}=-\frac{1}{\sin^2\theta}\mathcal{F}^{8}_{1214,24}, \\ 
	\mathcal{F}^{8}_{1434,13}=\frac{kr^4\sin^2\theta(c^2k+2U_2)}{4}=-\mathcal{F}^{8}_{1334,14},\\
	\mathcal{F}^{8}_{2434,23}=\frac{kr^4a^2\sin^2\theta(c^2k+2U_1)}{4c^2U}=-\mathcal{F}^{8}_{2334,24}; 
\end{cases}
		\end{eqnarray}

	\begin{eqnarray}\label{KK}
	\begin{cases}
		\mathcal{F}^{9}_{1223,13}=-\frac{kr^2(c^2k+2U_1)}{4U}=\frac{1}{\sin^2\theta}\mathcal{F}^{9}_{1224,14}=-\mathcal{F}^{9}_{1213,23}=\frac{1}{\sin^2\theta}\mathcal{F}^{9}_{1214,24}, \\ 
			\mathcal{F}^{9}_{1434,13}=\frac{kr^4\sin^2\theta(c^2k+2U_1)}{4}=-\mathcal{F}^{9}_{1334,14},\\
			\mathcal{F}^{9}_{2434,23}=\frac{kr^4a^2\sin^2\theta(c^2k+2U_1)}{4c^2U}=-\mathcal{F}^{9}_{2334,24};
		\end{cases}
				\end{eqnarray}


\begin{eqnarray}\label{Q(g,K)}
\begin{cases}	
	\mathcal{M}^{5}_{1223,13}=-\frac{3Y_6}{r^2}=\frac{1}{\sin^2\theta}\mathcal{M}^{5}_{1224,14}=-\mathcal{M}^{5}_{1213,23}=-\frac{1}{\sin^2\theta}\mathcal{M}^{5}_{1214,24},\\
	\mathcal{M}^{5}_{1434,13}=-\frac{Y_1Y_6\sin^2\theta}{r^2}=-\mathcal{M}^{5}_{1334,14}, \\
	\mathcal{M}^{5}_{2434,23}=-\frac{9r^2Y_6\sin^2\theta}{Y_1}=-\mathcal{M}^{5}_{2334,24}. \\
\end{cases}
		\end{eqnarray}

	
By examining the expressions in equations \eqref{RR} to
\eqref{Q(g,K)}, the following results are obtained:
	
	 \begin{pr}\label{pr5}
The FLRW spacetime satisfies the following pseudosymmetric type conditions:
			\begin{enumerate}[label=(\roman*)]
				\item
		   $\ C\cdot C=-\frac{k}{6a^2}Q(g,C),$
		
	\item  $C\cdot K=-\frac{k}{6a^2}Q(g,K),$
		
		\item  $K\cdot C=-\frac{c^2k+2U_1}{2c^2a^2}Q(g,C),$

\item  $K\cdot K=-\frac{c^2k+2U_1}{2c^2a^2}Q(g,K),$				
\item  $C\cdot K -K\cdot C=\frac{c^2k+3U_1}{3a^2c^2}Q(g,C),$
	\item  $R\cdot R-Q(S,R)=\frac{2a\ddot{a}}{c^2a}Q(g,C),$

\item  $C\cdot R+R\cdot C=-\frac{2c^2k+3U_1}{3c^2a^2}Q(g,C)+Q(S,C).$

		 	\end{enumerate}	
	\end{pr}

The non-zero components of the second-rank Ricci tensor $S^2_{qs}$ are derived as follows:
\begin{eqnarray}\label{S2}
\begin{cases}
S^2_{11}=\frac{9\ddot{a}^2}{c^2a},\ \ S^2_{22}=\frac{U_3^2}{c^4a^2U},\\
S^2_{33}=\frac{r^2(c^2k+U_3)^2}{c^4a^2}=\frac{1}{\sin^2\theta}S^2_{44}.\\
\end{cases}
\end{eqnarray}

Let $\mathcal{Y}^{1}=g\wedge g,$ $\mathcal{Y}^{2}= g\wedge S,$ $\mathcal{Y}^{3}=S\wedge S,$ $\mathcal{Y}^{4}=g\wedge S^2,$ $\mathcal{Y}^{5}=S\wedge S^2$ and $\mathcal{Y}^{6}=S^2\wedge S^2.$ Therefore the Kulkarni-Nomizu product tensors $\mathcal{Y}^{1},$ $\mathcal{Y}^{2},$ $\mathcal{Y}^{3},$ $\mathcal{Y}^{4},$ $\mathcal{Y}^{5},$ and $\mathcal{Y}^{6}$ are calculated as below:

\begin{eqnarray}\label{gag}
\begin{cases}	
 \mathcal{Y}^{1}_{1212}=-\frac{2c^2a^2}{U},\ \ \mathcal{Y}^{1}_{1313}=2c^2r^2a^2=\frac{1}{\sin^2\theta}\mathcal{Y}^{1}_{1414},\\
 \mathcal{Y}^{1}_{2323}=\frac{2r^2a^4}{U}=-\frac{1}{\sin^2\theta}\mathcal{Y}^{1}_{2424},\ \
	\mathcal{Y}^{1}_{3434}=2r^4a^4\sin^2\theta;\\
\end{cases}
		\end{eqnarray}

\begin{eqnarray}\label{Q(gas)}
\begin{cases}	
 \mathcal{Y}^{2}_{1212}=-\frac{2U_1}{U},\ \ \mathcal{Y}^{2}_{1313}=r^2(c^2k+2U_1)=\frac{1}{\sin^2\theta}\mathcal{Y}^{2}_{1414},\\
 \mathcal{Y}^{2}_{2323}=\frac{a^2r^2(c^2k+2U_1)}{c^2U}=\frac{1}{\sin^2\theta}\mathcal{Y}^{2}_{2424},\ \
	\mathcal{Y}^{2}_{3434}=-\frac{2r^4a^2\sin^2\theta(c^2k+U_3)}{c^2};\\
\end{cases}
		\end{eqnarray}
		
\begin{eqnarray}\label{Q(sas)}
\begin{cases}	
 \mathcal{Y}^{3}_{1212}=-\frac{6\ddot{a}U_3}{c^2aU},\ \ \mathcal{Y}^{3}_{1313}=\frac{6r^2\ddot{a}(c^2k+U_3)}{c^2a}=\frac{1}{\sin^2\theta}\mathcal{Y}^{3}_{1414},\\
 \mathcal{Y}^{3}_{2323}=\frac{2U_3r^2(c^2k+U_3)}{c^4U}=\frac{1}{\sin^2\theta}\mathcal{Y}^{3}_{2424},\ \
	\mathcal{Y}^{3}_{3434}=-\frac{2r^2\sin^2\theta(c^2k+U_3)}{c^4};\\
\end{cases}
		\end{eqnarray}
		
\begin{eqnarray}\label{Q(gas2)}
\begin{cases}	
 \mathcal{Y}^{4}_{1212}=-\frac{2(2\dot{a}^2U_1+5a^2\ddot{a})}{c^2a^2U},\ \ \mathcal{Y}^{4}_{1313}=\frac{9r^2a^2\ddot{a}^2+(c^2k+U_3)}{c^2a^2}=\frac{1}{\sin^2\theta}\mathcal{Y}^{4}_{1414},\\
 \mathcal{Y}^{4}_{2323}=\frac{r^2u_3^2+(c^2k+U_3)^2}{C^4U}=-\frac{1}{\sin^2\theta}\mathcal{Y}^{4}_{2424},\ \
	\mathcal{Y}^{4}_{3434}=-\frac{2r^4\sin^2\theta(c^2k+U_3)^2}{c^4};\\
\end{cases}
		\end{eqnarray}
		
\begin{eqnarray}\label{Q(sas2)}
\begin{cases}	
 \mathcal{Y}^{5}_{1212}=\frac{6\ddot{a}U_3U_4}{a^3c^4U},\ \ \mathcal{Y}^{5}_{1313}=\frac{3r^2\ddot{a}(c^2k+U_3)(c^2k+2U_4)}{c^4a^3}=\frac{1}{\sin^2\theta}\mathcal{Y}^{5}_{1414},\\
 \mathcal{Y}^{5}_{2323}=\frac{r^2U_3(c^2k+U_3)(c^2k+2U_4)}{c^6a^2U}=-\frac{1}{\sin^2\theta}\mathcal{Y}^{5}_{2424},\ \
	\mathcal{Y}^{5}_{3434}=\frac{2r^4\sin^2\theta(c^2k+U_3)^3}{c^6a^2};\\
\end{cases}
		\end{eqnarray}
		
\begin{eqnarray}\label{s2as2}
\begin{cases}	
 \mathcal{Y}^{6}_{1212}=-\frac{18\ddot{a}^2U_3}{a^2c^6U},\ \ \mathcal{Y}^{6}_{1313}=\frac{18\ddot{a}^2r^2(c^2k+U_3)^2}{a^4c^6}=\frac{1}{\sin^2\theta}\mathcal{Y}^{6}_{1414},\\
 \mathcal{Y}^{6}_{2323}=\frac{2r^2U_3^2(c^2k+U_3)^2}{a^4c^8U}=\frac{1}{\sin^2\theta}\mathcal{Y}^{6}_{2424},\ \
	\mathcal{Y}^{6}_{3434}=-\frac{2r^4\sin^2\theta(c^2k+U_3)^4}{a^4c^8}.\\
\end{cases}
		\end{eqnarray}	
Given $\eqref{gag}$ – $\eqref{s2as2}$, we have inferred the following linear dependency relation:

\begin{eqnarray}\label{GRT}
R&=&A_1 g\wedge g+A_2 g\wedge S+A_3 S\wedge S\nonumber\\ &+&A_4 g\wedge S^2+A_5 S\wedge S^2+A_6 S^2\wedge S^2,
\end{eqnarray}		
where $A_1,$ $A_2,$ $A_3,$ $A_4,$ $A_5,$ and  $A_6 $ are computed as follows:
$$A_1=1, A_2=1,$$	
\[
A_3=
\frac{
\begin{aligned}
& a\,\ddot a\Bigg[
- c^{2}\Big(
32 r^{4} \dot a^{4}
+ a\Big(3 c^{4} k r^{4}
+ 2 c^{2} (14 + 11 k r^{2}) r^{2}
+ 4 (3 + 2 c^{4} k r^{2}) r^{2} a^{2}\Big)
\dot a^{2} \big(c^{2} k r^{2} + 2 r^{2} \dot a^{2}\big) \\
& \qquad\qquad
+ 2 \big(c^{2} k r^{2} + 2 r^{2} \dot a^{2}\big)^{2}
\Big(c^{4} k r^{2} a^{2}
+ c^{2} (k r^{2} + 4 r^{2} a^{2}) \dot a^{2}
+ 2 r^{2} \dot a^{4}\Big)
\Big) \\
& \quad
+ r^{2} a \Big(
- 2 c^{4} k r^{2} (k r^{2} - 4 r^{2} a^{2})
+ c^{2} r^{2} (5 k r^{2} + 48 r^{2} a^{2}) \dot a^{2}
+ 36 r^{4} \dot a^{4}
\Big) \ddot a \\
& \quad
+ r^{4} a^{2} \Big(
c^{2} (23 k r^{2} + 32 r^{2} a^{2})
+ 94 r^{2} \dot a^{2}
\Big) \ddot a^{2}
+ 31 r^{6} a^{3} \ddot a^{3}
\Bigg]
\end{aligned}
}{
6\,\ddot a\,
\Big(
2\dot a^{2}
+ a\,\ddot a\, r^{6}
\big(c^{2}k+2U_{1}\big)
\big(c^{2}k+U_{3}\big)^{2}
\Big)
},
\]

\[
A_4=
\frac{
c^{2} r^{2} a^{2}\Big[
-2\big(c^{2} k r^{2} + 2 r^{2} \dot a^{2}\big)
\big(c^{2} k r^{2} + 2 r^{2} (c^{2} a^{2} + \dot a^{2})\big)
+ r^{2} a \big(c^{2} k r^{2} + 8 r^{2} (c^{2} a^{2} + \dot a^{2})\big)\ddot a
+ 9 r^{4} a^{2} \ddot a^{2}
\Big]
}{
2 r^{6}\big(c^{2}k+2U_{1}\big)\big(c^{2}k+U_{3}\big)^{2}
},
\]

\[
A_5=
\frac{
\begin{aligned}
& - c^{4} r^{2} a^{3}\Bigg[
-4\big(c^{2} k r^{2} + 2 r^{2} \dot a^{2}\big)^{2}
\Big(c^{4} k r^{2} a^{2}
+ c^{2} (k r^{2} + 4 r^{2} a^{2}) \dot a^{2}
+ 2 r^{2} \dot a^{4}\Big) \\
& \quad
- 2 a \big(c^{2} k r^{2} + 2 r^{2} \dot a^{2}\big)
\Big(c^{4} k r^{2} (3 k r^{2} + 8 r^{2} a^{2})
+ c^{2} r^{2} (19 k r^{2} + 24 r^{2} a^{2}) \dot a^{2}
+ 20 r^{4} \dot a^{4}\Big)\ddot a \\
& \quad
+ r^{2} a^{2}\Big(
- c^{4} k r^{4}
+ 2 r^{2}\big(8 c^{4} k r^{2} a^{2}
+ c^{2} (23 k r^{2} + 48 r^{2} a^{2}) \dot a^{2}
+ 72 r^{2} \dot a^{4}\big)
\Big)\ddot a^{2} \\
& \quad
+ r^{4} a^{3}\Big(
c^{2} (55 k r^{2} + 64 r^{2} a^{2})
+ 220 r^{2} \dot a^{2}
\Big)\ddot a^{3}
+ 68 r^{6} a^{4} \ddot a^{4}
\Bigg]
\end{aligned}
}{
6\ddot a\Big(
2\dot a^{2}
+ a\ddot a\, r^{8}
\big(c^{2}k+2U_{1}\big)
\big(c^{2}k+U_{3}\big)^{3}
\Big)
},
\]		

\[
A_6=
\frac{
\begin{aligned}
& - c^{6} r^{4} a^{5}\Bigg[
2\big(c^{2} k r^{2} + 2 r^{2} \dot a^{2}\big)
\Big(c^{4} k r^{2} a^{2}
+ c^{2} (k r^{2} + 4 r^{2} a^{2}) \dot a^{2}
+ 2 r^{2} \dot a^{4}\Big) \\
& \quad
+ a\Big(
3 c^{4} k r^{2} (k r^{2} + 2 r^{2} a^{2})
+ 4 c^{2} r^{2} (2 k r^{2} + r^{2} a^{2}) \dot a^{2}
- 2 r^{4} \dot a^{4}
\Big)\ddot a \\
& \quad
- r^{2} a^{2}\Big(
c^{2} (7 k r^{2} + 20 r^{2} a^{2})
+ 38 r^{2} \dot a^{2}
\Big)\ddot a^{2}
- 22 r^{4} a^{3} \ddot a^{3}
\Bigg]
\end{aligned}
}{
6\ddot a\Big(
2\dot a^{2}
+ a\ddot a\, r^{8}
\big(c^{2}k+2U_{1}\big)
\big(c^{2}k+U_{3}\big)^{3}
\Big)
}.
\]		

Contracting the last relation yields the following expression:
\begin{equation}\label{Ein3}
S^3+\frac{\dot{a}^2-k-5U_1}{c^2r^2a^2}S^2+\frac{-3a(c^2k+2U_3)}{c^6r^2a^5}S+\frac{3\ddot{a}U_3(c^2k+U_3)}{c^6r^2a^5}g=0.
\end{equation}
This gives rise to the following:
\begin{pr}\label{pr2}
\begin{enumerate}[label=(\roman*)]

The nature of the FLRW spacetime is neither an Einstein manifold of level $2$ nor Roter type, but rather it is 
\item an Einstein manifold of level $3,$
\item generalized Roter type as shown in \eqref{GRT},
\item $2$- quasi Eintein as the rank of $(S-\zeta g)$ is $2$ for $\zeta=\frac{c^2k+U_3}{c^2a^2}.$
\end{enumerate}
\end{pr}	

\begin{pr}\label{Com5}
The FLRW spacetime admits the following curvature conditions:
\begin{enumerate}
\item The general form of a tensor compatible with $R$ is
$$
		\left(
		\begin{array}{cccc}
			\mathscr{Z}_{11} &\mathscr{Z}_{12} & 0 & 0 \\
			\frac{\dot{a}^2}{a\ddot{a}}\mathscr{Z}_{12} & \mathscr{Z}_{22} & 0 & 0 \\
			0 & 0 & \mathscr{Z}_{33} & \mathscr{Z}_{34} \\
			0 & 0 & \mathscr{Z}_{34} & \mathscr{Z}_{44}
		\end{array}
		\right),
		$$
\item The common form of conformal $(C)$ and conharmonic $(K)$ compatible tensor are
$$
		\left(
		\begin{array}{cccc}
			\mathscr{Z}_{11} &\mathscr{Z}_{12} & 0 & 0 \\
			\mathscr{Z}_{12} & \mathscr{Z}_{22} & 0 & 0 \\
			0 & 0 & \mathscr{Z}_{33} & \mathscr{Z}_{34} \\
			0 & 0 & \mathscr{Z}_{34} & \mathscr{Z}_{44}
		\end{array}
		\right),
		$$		
\item The general form for a tensor compatible with Projective $(P)$ is
$$
		\left(
		\begin{array}{cccc}
			\mathscr{Z}_{11} &\frac{c^2k+2U_1}{c^2k-U_1}\mathscr{Z}_{21} & 0 & 0 \\
			\mathscr{Z}_{21} & \mathscr{Z}_{22} & \mathscr{Z}_{23} & \mathscr{Z}_{24} \\
			0 & \frac{-2(c^2k+U_2)}{c^2k-U_2}\mathscr{Z}_{23} & \mathscr{Z}_{33} & \mathscr{Z}_{34} \\
			0 & \frac{-2(c^2k+U_2)}{c^2k-U_2}\mathscr{Z}_{24} & \mathscr{Z}_{34} & \mathscr{Z}_{44}
		\end{array}
		\right),
		$$
\item The common form of a concircular $(W)$ compatible tensor is
$$
		\left(
		\begin{array}{cccc}
			\mathscr{Z}_{11} &\mathscr{Z}_{12} & 0 & 0 \\
			\frac{c^2k-3U_2}{c^2k+3U_2}\mathscr{Z}_{12} & \mathscr{Z}_{22} & 0 & 0 \\
			0 & 0 & \mathscr{Z}_{33} & \mathscr{Z}_{34} \\
			0 & 0 & \mathscr{Z}_{34} & \mathscr{Z}_{44}
		\end{array}
		\right),
		$$				
\end{enumerate}
\end{pr}
\begin{thm}
The following geometric characteristics are possessed by the FLRW spacetime:

	\begin{enumerate}[label=(\roman*)]
		
		\item it is conformal pseudosymmetric due to conformal curvature tensor as  $\ C\cdot C=-\frac{k}{6a^2}Q(g,C),$ \\
		\item it is conharmonic pseudosymmetric due to conformal curvature as  $C\cdot K=-\frac{k}{6a^2}Q(g,K),$ \\
%
%
%
		 
	\item it exhibits conformal pseudosymmetric due to conharmonic curvature as $K\cdot C=-\frac{c^2k+2U_1}{2c^2a^2}Q(g,C);$
\item it realizes the pseudosymmetric type curvature relation $C\cdot K -K\cdot C=\mathscr{L}_1Q(g,C),$	where $\mathscr{L}_1=	\frac{c^2k+3U_1}{3a^2c^2},$		
		  
\item it satisfies the pseudosymmetric type curvature condition $R\cdot R-Q(S, R)=\mathscr{L}_2Q(g,C),$ where $\mathscr{L}_2=\frac{2a\ddot{a}}{c^2a},$
		
	\item the tensors $Q(g, C)$ and $Q(S, C),$ can be linearly combined to form the tensor $C\cdot R+R\cdot C,$

		 

	\item for $\zeta=\frac{c^2k+U_3}{c^2a^2},$ it is a $2$-quasi Einstein manifold,

		\item it is an Einstein manifold of level $3,$ (see, $\eqref{Ein3}),$
		  
\item it is generalized Roter type spacetime $(see, \eqref{GRT}),$ 
		
		
		 
		\item the universal form of tensors in FLRW spacetime that are compatible with $R$, $C$, $P$, $W$ and $K$ are derived (see, \eqref{Com5}).

	\end{enumerate}
	
\end{thm}

\begin{rem}
We observe that FLRW spacetime does not adhere to any of the following by examining different components of various curvature conditions \eqref{RR} to \eqref{Q(g,K)}.
	\begin{enumerate}[label=(\roman*)]
		\item considering that $\nabla R\neq 0$, it satisfies the conditions $\nabla W\neq 0$, $\nabla P\neq 0$, $\nabla K\neq 0$, $\nabla C\neq 0$.
		


		
		
		
		\item for any curvatures $R, W, P, C$, and $K$, it is not Venzi-space. 
		
\item 	the Ricci tensor of the FLRW spacetime is neither cyclic parallel nor Codazzi-type,	
		
		\item there is neither weak symmetry nor Chaki pseudosymmetry in the FLRW spacetime for $R$ $W,$ $P$ $C,$ and $K$.
		
	\end{enumerate}  
\end{rem}

\section{\bf Ricci soliton on Friedmann–Lemaître–Robertson–Walker spacetime}

 In the FLRW spacetime, the vector field $\frac{\partial}{\partial \phi}$ acts as a Killing vector field, which implies that the Lie derivative of the metric tensor satisfies $\pounds_{\partial/\partial \phi} g = 0,$ where $\pounds$ denotes the Lie derivative operator. On the other hand, the coordinate vector fields $\frac{\partial}{\partial t}$, $\frac{\partial}{\partial r}$, and $\frac{\partial}{\partial \theta}$ are non-Killing. 
 
 Let us define $\mathcal{A}^1 = \pounds_{\frac{\partial}{\partial t}} g$ and $\mathcal{A}^2 = \pounds_{\frac{\partial}{\partial r}} g$. The non-vanishing components of the tensors $\mathcal{A}^1$ and $\mathcal{A}^2$ are obtained as follows:
$$\begin{array}{c}
	\mathcal{A}^1_{22}=\frac{2a\dot{a}}{U},\ \ \ 
	\mathcal{A}^1_{33}=-2r^2a\dot{a},\ \ \
\mathcal{A}^1_{44}=2r^2a\dot{a}\sin^2\theta.\nonumber\\
	\mathcal{A}^2_{33}=-2a^2,\ \ \
		\mathcal{A}^2_{44}=-2ra^2\sin^2\theta.\nonumber
	\nonumber\\
\end{array}$$

Hence, in the FLRW spacetime, for the non-Killing vector field $\frac{\partial}{\partial t}$ and the 1-form $\eta = \left(0,\frac{ck\sqrt{2a\dot{a}}}{\sqrt{c^2k(1+r^2)-2UU_1}},0,0\right)$, the following relation is satisfies:

\begin{eqnarray}\nonumber
	\pounds_{\frac{\partial}{\partial t}}g+ 2\sigma_1 S+2\sigma_2 g+2\sigma_3 (\eta\otimes\eta)=0,
\end{eqnarray}
where $\sigma_1,\sigma_2,\sigma_3$ are given by
\begin{equation}\label{Coefficient1}
	\left.
	\begin{aligned}
		\sigma_1&=-\frac{c^2a\dot{a}}{c^2k+2U_2},\\
		\sigma_2&=-\frac{6a\dot{a}}{c^2k+2U_2},\\
		\sigma_3&=-1.
	\end{aligned}\ \ 
	\right\rbrace	
\end{equation}

Also, in the FLRW spacetime, considering the non-Killing vector field $\frac{\partial}{\partial r}$ and the $1$-form $\eta=\left(\frac{2\sqrt{r(a\ddot{a}-\dot{a}^2)}}{kr},0,0,0\right)$, the following relation is satisfies:

\begin{eqnarray}\nonumber
	\pounds_{\frac{\partial}{\partial r}}g+ 2\sigma_4 S+2\sigma_5 g+2\sigma_6 (\eta\otimes\eta)=0,
\end{eqnarray}
where $\sigma_4,\sigma_5,\sigma_6$ are given by
\begin{equation}\label{Coefficient2}
	\left.
	\begin{aligned}
		\sigma_4&=-\frac{2a^2}{kr},\\
		\sigma_5&=-\frac{U_3}{c^2kr},\\
		\sigma_6&=-1.
	\end{aligned}\ \ 
	\right\rbrace	
\end{equation}
This leads to the following:

\begin{thm} \label{thm6.1}
	In the FLRW spacetime, if the quantity $(c^2k + 2(\dot{a}^2 - a\ddot{a}))$ is non-vanishing, then the spacetime admits an almost $\eta$-Ricci–Yamabe soliton with soliton vector field $\xi = \frac{\partial}{\partial t}$ and the associated $1$-form $\eta = \left(0,\frac{ck\sqrt{2a\dot{a}}}{\sqrt{c^2k(1+r^2)-2UU_1}},0,0\right)$ satisfying:
	$$\frac{1}{2}\pounds_\xi g+\sigma_1 S+\left(\lambda-\frac{1}{2}\sigma_7\kappa\right) g+\sigma_3 (\eta\otimes\eta)=0,$$
	where $\sigma_7=2$,   $\lambda=\sigma_2+\kappa$, and $\sigma_1,\sigma_2,\sigma_3$ are given in (\ref{Coefficient1}).
	
\end{thm}

\begin{thm}\label{thm6.2}
In the \textit{FLRW spacetime}, if
$c^2\big(k+a\ddot{a}+2(\dot{a}^2-a\ddot{a})\big)=0
\quad \text{and} \quad 
(c^2k+2(\dot{a}^2-a\ddot{a})) \neq 0,$
then the spacetime admits an almost $\eta$-Ricci soliton with soliton vector field $\xi=\frac{\partial}{\partial t}$ and the associated $1$-form
$\eta=\left(0,\frac{ck\sqrt{2a\dot{a}}}{\sqrt{c^2k(1+r^2)-2UU_1}},0,0\right)$
satisfying
	\begin{eqnarray}
		&&\frac{1}{2}\pounds_{\xi}g+S+\sigma_2 g-\sigma_3(\eta\otimes\eta)=0\nonumber,
	\end{eqnarray}
where $\sigma_2$ and $\sigma_3$ are given in (\ref{Coefficient1}).
\end{thm}

\begin{thm} \label{thm6.3}
	In the FLRW spacetime, if the quantity $c^2kr$ is non-zero, then the spacetime admits an almost $\eta$-Ricci--Yamabe soliton with soliton vector field $\xi = \frac{\partial}{\partial r}$ and the associated $1$-form $\eta = \left(\frac{2\sqrt{r(a\ddot{a}-\dot{a}^2)}}{kr},0,0,0\right),$
	which satisfies the following relation:
	$$\frac{1}{2}\pounds_\xi g+\sigma_4 S+\left(\lambda-\frac{1}{2}\sigma_8\kappa\right) g+\sigma_6 (\eta\otimes\eta)=0,$$
	where $\sigma_8=2$,   $\lambda=\sigma_5+\kappa$, and $\sigma_4,\sigma_5,\sigma_6$ are given in (\ref{Coefficient2}).
	
\end{thm}

\begin{thm}\label{thm6.4}
In the \textit{FLRW spacetime}, assume that
$kr+2a^2=0 \quad \text{and} \quad c^2kr \neq 0.$
Then the spacetime admits an almost $\eta$-Ricci soliton with soliton vector field $\xi=\frac{\partial}{\partial r}$ and the associated $1$-form
$
\eta=\left(\frac{2\sqrt{r(a\ddot{a}-\dot{a}^2)}}{kr},0,0,0\right)
$
satisfying
	\begin{eqnarray}
		&&\frac{1}{2}\pounds_{\xi}g+S+\sigma_5 g-\sigma_6(\eta\otimes\eta)=0\nonumber,
	\end{eqnarray}
where $\sigma_5$ and $\sigma_6$ are given in (\ref{Coefficient2}).
\end{thm}

\begin{rem}
From Theorem \ref{thm6.2}, for the soliton vector field $\frac{\partial}{\partial t},$ the following cases arise for the FLRW spacetime:\\
\textbf{Case 1:} If $k=0$, then the FLRW spacetime admits an \emph{almost Ricci soliton}.\\
\textbf{Case 2:} If $a(t)$ is a non-zero constant, then the FLRW spacetime admits a \emph{Ricci soliton}.

\end{rem}

\begin{rem}
It follows from Theorem \ref{thm6.4} that if the FLRW spacetime satisfies the condition 
$a(t)=C_1 e^{C_2t}$, $C_1, C_2$ being constants, then the spacetime admits an almost Ricci soliton with respect to the soliton vector field $\frac{\partial}{\partial r}$.
\end{rem}
\section{\bf Generalized curvature inheritance on Friedmann–Lemaître–Robertson–Walker spacetime}
For a $(1,3)$-type Riemann curvature tensor, the concept of curvature collineation was first introduced by Katzin \textit{et al.} in 1969 \cite{KLD1969, KLD1970}, and it is characterized by the condition $\pounds_\xi \widetilde{R}=0$, where $\widetilde{R}$ denotes the $(1,3)$-type Riemann curvature tensor. Subsequently, Duggal \cite{Duggal1992} generalized this concept in 1992 by formulating the notion of curvature inheritance, which is the relation $\pounds_\xi \widetilde{R}=\gamma \widetilde{R}$. During the past three decades, a considerable number of consequences (see, \cite{Ahsan1978, Ahsan1977_231, Ahsan1977_1055, Ahsan1987, Ahsan1995, Ahsan1996, Ahasan2005, AhsanAli2014, AA2012, AH1980, AliAhsan2012, SASZ2022, ShaikhDatta2022}) have been investigated and expanded on different aspects of these symmetry properties.

Let $\pounds_{\frac{\partial}{\partial t}}R=I,$ $\pounds_{\frac{\partial}{\partial r}}R=I^1,$ the non-zero components of $I,$ $I^1$ as follows:
\begin{eqnarray}\label{LvR}
\begin{cases}	
I_{1212}=-\frac{\dot{a}\ddot{a}+aa^{(3)}}{U}=-I_{1221}=-I_{2112}=I_{2121}, \\
I_{1313}=r^2(\dot{a}\ddot{a}+aa^{(3)})=-I_{1331}=\frac{1}{\sin^2\theta}I_{1414}=-\frac{1}{\sin^2\theta}I_{1441}=I_{3113}=I_{3131}\\=-\frac{1}{\sin^2\theta}I_{4114}=\frac{1}{\sin^2\theta}I_{4141}, \\
I_{2323}=\frac{2r^2a\dot{a}U_1}{c^2U}=I_{2332}=\frac{1}{\sin^2\theta}I_{2424}=-\frac{1}{\sin^2\theta}I_{2442}=I_{3223}=I_{3232}=-\frac{1}{\sin^2\theta}I_{4224}=\frac{1}{\sin^2\theta}I_{4242}, \\
I_{3434}=-\frac{2r^4a\dot{a}\sin^2\theta(c^2k+U_1)}{c^2}=I_{3443}=\frac{1}{\sin^2\theta}I_{4334}=-\frac{1}{\sin^2\theta}I_{4343}, \\
\end{cases}
		\end{eqnarray}
\begin{eqnarray}\label{LtR}
\begin{cases}	
I^1_{1313}=2ra\ddot{a}=-I^1_{1331}=\frac{1}{\sin^2\theta}I^1_{1414}=-\frac{1}{\sin^2\theta}I^1_{1441}=I^1_{3113}=I^1_{3131}=-\frac{1}{\sin^2\theta}I^1_{4114}=\frac{1}{\sin^2\theta}I^1_{4141}, \\
I^1_{2323}=\frac{2ra^2\dot{a}^2}{c^2U}=-I^1_{2332}=\frac{1}{\sin^2\theta}I^1_{2424}=-\frac{1}{\sin^2\theta}I^1_{2442}=-I^1_{3223}=I^1_{3232}=-\frac{1}{\sin^2\theta}I^1_{4224}=\frac{1}{\sin^2\theta}I^1_{4242}, \\
I^1_{3434}=-\frac{2r^3a^2\sin^2\theta(c^2k+2\dot{a}^2)}{c^2}=I^1_{3443}=-I^1_{4334}=I^1_{4343}. \\
\end{cases}
		\end{eqnarray}

If $V=\nabla t =\{\frac{1}{c^2},0,0,0\},$ and $\pounds_ {V}R=I^2,$ $\pounds_{V}C=I^3,$ $\pounds_{V}W=I^4,$ $\pounds_{V}K=I^5,$ then the non-zero components of $I^2,$ $I^3,$ $I^4$  and $I^5$ are as follows:
\begin{eqnarray}\label{LvR}
\begin{cases}	
I^2_{1212}=-\frac{\dot{a}\ddot{a}+aa^{(3)}}{c^2U}=-I^2_{1221}=-I^2_{2112}=I^2_{2121}, \\
I^2_{1313}=\frac{r^2(\dot{a}\ddot{a}+aa^{(3)})}{c^2}=-I^2_{1331}=\frac{1}{\sin^2\theta}I^2_{1414}=-\frac{1}{\sin^2\theta}I^2_{1441}=I^2_{3113}=I^2_{3131}=-\frac{1}{\sin^2\theta}I^2_{4114}=\frac{1}{\sin^2\theta}I^2_{4141}, \\
I^2_{2323}=\frac{2r^2a\dot{a}U_1}{c^4U}=I^2_{2332}=\frac{1}{\sin^2\theta}I^2_{2424}=-\frac{1}{\sin^2\theta}I^2_{2442}=I^2_{3223}=I^2_{3232}=-\frac{1}{\sin^2\theta}I^2_{4224}=\frac{1}{\sin^2\theta}I^2_{4242}, \\
I^2_{3434}=-\frac{2r^4a\dot{a}\sin^2\theta(c^2k+U_1)}{c^4}=I^2_{3443}=\frac{1}{\sin^2\theta}I^2_{4334}=-\frac{1}{\sin^2\theta}I^2_{4343}, \\
\end{cases}
		\end{eqnarray}
	
\begin{eqnarray}\label{LvC}
\begin{cases}	
I^3_{2323}=-\frac{kr^2a\dot{a}}{3c^2U}=I^3_{2332}=\frac{1}{\sin^2\theta}I^3_{2424}=-\frac{1}{\sin^2\theta}I^3_{2442}=I^3_{3223}=I^3_{3232}=-\frac{1}{\sin^2\theta}I^3_{4224}=\frac{1}{\sin^2\theta}I^3_{4242}, \\
I^3_{3434}=-\frac{2r^4ka\dot{a}\sin^2\theta}{3c^2}=I^3_{3443}=\frac{1}{\sin^2\theta}I^3_{4334}=-\frac{1}{\sin^2\theta}I^3_{4343}, \\
\end{cases}
		\end{eqnarray}
		
\begin{eqnarray}\label{LvW}
\begin{cases}	
I^4_{1212}=\frac{\dot{a}\ddot{a}-aa^{(3)}}{2c^2U}=-I^4_{1221}=-I^4_{2112}=I^4_{2121}, \\
I^4_{1313}=-\frac{r^2(\dot{a}\ddot{a}-aa^{(3)})}{2c^2}=-I^4_{1331}=\frac{1}{\sin^2\theta}I^4_{1414}=-\frac{1}{\sin^2\theta}I^4_{1441}=I^4_{3113}=I^4_{3131}\\=-\frac{1}{\sin^2\theta}I^4_{4114}=\frac{1}{\sin^2\theta}I^4_{4141}, \\
I^4_{2323}=-\frac{ar^2(6\dot{a}^3+\dot{a}(c^2k+3a\ddot{a})+3a^2aq^{(3)})}{6c^4U}=I^4_{2332}=\frac{1}{\sin^2\theta}I^4_{2424}=-\frac{1}{\sin^2\theta}I^4_{2442}=I^4_{3223}=I^4_{3232}\\
=-\frac{1}{\sin^2\theta}I^4_{4224}=\frac{1}{\sin^2\theta}I^4_{4242}, \\
I^4_{3434}=\frac{r^4a\sin^2\theta(\dot{a}(-10c^2k-3r^2(a\ddot{a}-2\dot{a}^2)+3a^2a^{(3)}))}{6c^4}=I^4_{3443}=\frac{1}{\sin^2\theta}I^4_{4334}=-\frac{1}{\sin^2\theta}I^4_{4343}, \\
\end{cases}
		\end{eqnarray}	
		
\begin{eqnarray}\label{LvK}
\begin{cases}	
I^5_{1212}=\frac{3\dot{a}\ddot{a}+aa^{(3)}}{c^2U}=-I^5_{1221}=-I^5_{2112}=I^5_{2121}, \\
I^5_{1313}=-\frac{r^2(3\dot{a}\ddot{a}+aa^{(3)})}{c^2}=-I^5_{1331}=\frac{1}{\sin^2\theta}I^5_{1414}=-\frac{1}{\sin^2\theta}I^5_{1441}=I^5_{3113}=I^5_{3131}\\=-\frac{1}{\sin^2\theta}I^5_{4114}=\frac{1}{\sin^2\theta}I^5_{4141}, \\
I^5_{2323}=-\frac{a\dot{a}r^2(2\dot{a}^2+(c^2k+5a\ddot{a})+a^2a^{(3)})}{c^4U}=\frac{1}{\sin^2\theta}I^5_{2424}=-\frac{1}{\sin^2\theta}I^5_{2442}=I^5_{3232}=-\frac{1}{\sin^2\theta}I^5_{4224}=\frac{1}{\sin^2\theta}I^5_{4242}, \\
I^5_{2332}=\frac{ar^2(\dot{a}(c^2k+(2\dot{a}^2+5a\ddot{a}))+a^3a^{(3)})}{c^4U}=I^5_{3223},\\
I^5_{3434}=\frac{r^4a\sin^2\theta(2\dot{a}^3+5a\dot{{a}\ddot{a}+a^2a^{(3)}})}{c^4}=I^5_{3443}=\frac{1}{\sin^2\theta}I^5_{4334}=-\frac{1}{\sin^2\theta}I^5_{4343}. \\
\end{cases}
		\end{eqnarray}		

\begin{thm}
The FLRW spacetime realizes the following inheritance properties
\begin{enumerate}[label=(\roman*)]
\item Generalized curvature inheritance, i.e., $\pounds_{\frac{\partial}{\partial t}}R=p_1R+p_2 g\wedge g+p_3g\wedge S+p_4S\wedge S,$ \\
where $$p_1=\frac{2\dot{a}^3-5a\dot{a}\ddot{a}+a^2a^{(3)}}{aU_2},$$
\begin{equation*}
p_2=\frac{
\begin{aligned}
&-2 \dot a^5\, \ddot a
+ \dot a^3 \ddot a \left(-c^2 k + 10 a \ddot a \right)
+ a \dot a \ddot a^2 \left(4 c^2 k + 19 a \ddot a \right) \\
&+ 2 a \dot a^4 a^{(3)}
+ a \dot a^2 \left(c^2 k - 6 a \ddot a \right) a^{(3)}
- a^2 \ddot a \left(2 c^2 k + 5 a \ddot a \right) a^{(3)}
\end{aligned}
}{
2 c^2 a U_2 (c^2 k + 2 U_2)
},
\end{equation*}
$$p_3=-\frac{3a\ddot{a}(3\dot{a}\ddot{a}-aa^{(3)})}{U_2 (c^2 k + 2 U_2)},$$
$$p_4=\frac{3a^2c(3\dot{a}\ddot{a}+aa^{(3)})}{2U_2 (c^2 k + 2 U_2)}.$$
\item Generalized curvature inheritance, i.e., $\pounds_{\frac{\partial}{\partial r}}R=p_5R+p_6 g\wedge g+p_7g\wedge S+p_8S\wedge S,$  \\
where $$p_5=\frac{2(c^2k+2\dot{a}^2)}{c^kr},$$
$$p_6=\frac{\ddot a \left(
c^4 k
+ 8 \dot a^4
+ a \ddot a \left(2 c^2 k + a \ddot a \right)
+ \dot a^2 \left(6 c^2 k + 9 a \ddot a \right)
\right)}{c^4rk(c^2k+2U_2)},$$
$$p_7=\frac{2a\dot{a}(c^2k+5\dot{a}^2+a\ddot{a})}{c^2kr(c^2k+2U_2)},$$
$$p_8=-\frac{a^2U}{k^2r(c^2k+2U_2)}.$$
\end{enumerate}
\end{thm}			
\begin{thm}\label{GCIG}
For the non-Killing vector field $V=\nabla t=\frac{1}{c^2}\frac{\partial}{\partial t},$ the FLRW spacetime explores the following properties:
\begin{enumerate}[label=(\roman*)]
\item Generalized curvature inheritance, i.e., $\pounds_{V}R=\mu_1R+\mu_2 g\wedge g+\mu_3g\wedge S+\mu_4S\wedge S,$\\
where
 $$\mu_1=\frac{2\dot{a}^2-5a\dot{a}\ddot{a}+a^2a^{(3)}}{c^2aU_1},$$

\begin{equation*}
\mu_2=\frac{
\begin{aligned}
&-2 \dot a^5\, \ddot a
+ \dot a^3 \ddot a \left(-c^2 k + 10 a \ddot a \right)
+ a \dot a \ddot a^2 \left(4 c^2 k + 19 a \ddot a \right) \\
&+ 2 a \dot a^4 a^{(3)}
- a^2 \ddot a \left(2 c^2 k + 5 a \ddot a \right) a^{(3)}
+ \dot a^2 \left(c^2 k\, a - 6 a^2 \ddot a \right) a^{(3)}
\end{aligned}
}{
2c^4 a^2 U_2 (c^2 k + 2 U_2)
},
\end{equation*}

$$\mu_3=-\frac{3a\ddot{a}(3\dot{a}\ddot{a}-aa^{(3)})}{c^2U_2(c^2k+2U_2)},$$
$$\mu_4=\frac{a^2(3\dot{a}\ddot{a}-aa^{(3)})}{2U_2(c^2k+\dot{a}^2-2a\ddot{a})}.$$

\item Generalized conformal curvature inheritance, i.e., $\pounds_{X}C=\mu_5C+\mu_6 g\wedge g+\mu_7g\wedge S+\mu_8S\wedge S,$  
where $$\mu_5=\frac{\dot{a}(c^2k+6U_2)}{3c^2aU_2},$$
$$\mu_6=\frac{\dot{a}r^2(c^4(c^2k+11k\dot{a}^2-2a\ddot{a})+c^2(36k\dot{a}^4-9ka^2\ddot{a}^2)+36\dot{a}^6-27a^2\dot{a}^2\ddot{a}^2-9a^3\ddot{a}^3)}{18c^4\dot{a}^2U_2(c^2k+2U_2)},$$
$$\mu_7=\frac{\dot{a}(c^4k+8c^2k\dot{a}^2+12\dot{a}^4-2c^2ka\ddot{a}-6a\dot{a}^2\ddot{a}-6a^2\ddot{a}62)}{6c^2aU_2(c^2k+2U_2)},$$
$$\mu_8=\frac{a\dot{a}(c^2k+3U_2)}{6U_2(c^2k+2U_2)}.$$
\item Generalized concicular curvature inheritance, i.e., $\pounds_{X}W=\mu_9W+\mu_{10} g\wedge g+\mu_{11}g\wedge S+\mu_{12}S\wedge S,$
where
$$\mu_9=\frac{c^2k\dot{a}+9\dot{a}^3-12a\dot{a}\ddot{a}+3a^2a^{(3)}}{3c^2aU_2},$$
\begin{equation*}
\mu_{10}=\frac{
\begin{aligned}
&54 \dot a^7
+ 3 \dot a^5 \left(17 c^2 k - 66 a \ddot a \right) 
+ 2 \dot a^3 \left(7 c^4 k
+ 9 a \ddot a \left(-7 c^2 k + 4 a \ddot a \right)\right) \\
&+ \dot a \left(c^6 k
+ a \ddot a \left(-17 c^4 k
+ 3 a \ddot a \left(7 c^2 k + 78 a \ddot a \right)\right)\right)
+ 36 a^2 \dot a^4 a^{(3)}
+ 24 a^2 \dot a^2 \left(c^2 k - 3 a \ddot a \right) a^{(3)} \\
&+ 3 a^2 \left(c^4 k
- 2 a \ddot a \left(4 c^2 k + 21 a \ddot a \right)\right) a^{(3)}
\end{aligned}
}{
36 c^4 a^3 U_2 (c^2 k + 2 U_2)
},
\end{equation*}
 
$$\mu_{12}=\frac{c^2ka\dot{a}+3a\dot{a}^3-6a^2\dot{a}^2\ddot{a}+3a^3a^{(3)}}{6U_2(c^2 k + 2 U_2}.$$
\item Generalized conharmonic curvature inheritance,  i.e., $\pounds_{X}K=\mu_{13}K+\mu_{14} g\wedge g+\mu_{15}g\wedge S+\mu_{16}S\wedge S,$
where $$\mu_{13}=\frac{\dot{a}(c^2k+4\dot{a}^2)}{c^2aU_2},$$
$$\mu_{14}=\frac{c^4k\dot{a}+6c^2k\dot{a}^3+8\dot{a}^5+4c^2ka\dot{a}\ddot{a}+10\dot{a}^3\ddot{a}+6a^2\dot{a}\ddot{a}^2}{2c^2aU_2(c^2k+2U_2)},$$
$$\mu_{15}=-\frac{c^2ka\dot{a}+3a\dot{a}^3+a^2\dot{a}\ddot{a}}{2U_2(c^2k+2U_2)},$$
$$\mu_{16}=\frac{\dot{a}\ddot{a}(c^4kr^2+(14\dot{a}^4+\dot{a}^2(9c^2k+9a\ddot{a}(-7c^2k+4a\ddot{a})))+r^2aU_2(c^2k+2U_2))a^{(3)}}{2c^4a^2U_2(c^2k+2U_2)},$$

and $a$ represents a function of $t$.
\end{enumerate}
\end{thm}	
\begin{cor}
If $a(t)$ is a constant, then the Theorem \ref{GCIG} does not hold.
\end{cor}
The non-zero components of $\pounds_{\frac{\partial}{\partial \theta}}S$  and $\pounds_{\frac{\partial}{\partial \theta}}g$ are calculated as follows:\\
$\left(\pounds_{\frac{\partial}{\partial \theta}}S\right)_{44}=-\frac{r^2\sin 2\theta(c^2k+U_3)}{c^2},$\\
$\left(\pounds_{\frac{\partial}{\partial \theta}}g\right)_{44}=-r^2a^2\sin 2\theta$.
\begin{thm}
For the non-Killing soliton vector field $\frac{\partial}{\partial \theta}$, the FLRW spacetime realizes the inheritance symmetry.\\
i.e.,
$$\pounds_{\frac{\partial}{\partial \theta}}g=\frac{c^2a^2}{c^2k+U U3}\pounds_{\frac{\partial}{\partial \theta}}S.$$ 
\end{thm}
\begin{rem}
The FLRW spacetime does not satisfies the following structures:
\begin{enumerate}[label=(\roman*)]
\item it does not fulfill any curvature collineation for the non-Killing vector field $\frac{\partial}{\partial t},$ $\frac{\partial}{\partial r}$ and $\frac{\partial}{\partial \theta}.$
\item it does not satisfies Ricci inheritance symmetries.
\end{enumerate}
\end{rem}
\section{\bf Friedmann–Lemaître–Robertson–Walker spacetime Vs Lemaître–Tolman–Bondi spacetime  }
The Lemaître–Tolman–Bondi (briefly, LTB) spacetime \cite{Bondi1947, Lem1997, Tolman1934} serves as a simplified model for representing an inhomogeneous universe. One of the main advantages of the LTB metric is that it provides a spherically symmetric exact solution of the Einstein field equations. Moreover, this model has been considered as a possible alternative to the concept of dark energy within the framework of the standard cosmological model. In a spherical coordinate system, the corresponding metric can be written as follows (see, \cite{SAACD_LTB_2022}, Eqn $(1.1)$):
\begin{eqnarray}\label{VBdS}
	ds^2=-dt^2+\frac{B^{\prime^2}(t,r)}{1+2\mu(r)}dr^2+B^2(t,r)\left(d\theta^2+\sin^2\theta d\phi^2\right),
    \end{eqnarray}
     where $B(t, r )$ is the areal radius function, $\mu(r )$ is the curvature function and $B^{\prime}(t,r)=\frac{\partial}{\partial r} B(t,r).$  A comprehensive analysis of the Lemaître–Tolman–Bondi spacetime is presented in \cite{SAACD_LTB_2022}. Here, a comparative overview of the curvature and physical characteristics of FLRW spacetime and the LTB spacetime is summarized below:\\
%

\begin{table}[H]
\centering
\caption{\textbf{Similarities:}}
\label{tab:flrw_ltb_similarities}
\renewcommand{\arraystretch}{1.5}
\begin{tabular}{|p{0.3\textwidth}|p{0.65\textwidth}|}
\hline
\textbf{Property} & \textbf{Description of similarity} \\
\hline
\textbf{Scalar curvature} & The scalar curvature takes a non-vanishing value in both spacetimes. \\
\hline
\textbf{Manifold classification} & Both spacetimes belong to the class of $Ein(3)$ manifolds. \\
\hline
\textbf{Quasi-Einstein specification} & Both spacetimes specify $2$-quasi-Einstein manifolds. \\
\hline
\textbf{Killing vector fields} & For both spacetimes, the coordinate vector fields $\frac{\partial}{\partial \phi}$ are Killing. \\
\hline
\textbf{Non-Killing vector fields} & In both spacetimes, the vector fields $\frac{\partial}{\partial t}$, $\frac{\partial}{\partial r}$, and $\frac{\partial}{\partial \theta}$ are non-Killing. \\
\hline
\textbf{Ricci tensor classification} & The Ricci tensors of both spacetimes are neither of the cyclic parallel nor Codazzi type. \\
\hline
\textbf{Tensor compatibility} & For both spacetimes, the Ricci tensor is conharmonic, Riemann, and concircular compatible. \\
\hline
\textbf{Almost Ricci solitons} & For the non-Killing vector field $\frac{\partial}{\partial r}$, both spacetimes admit almost Ricci soliton under particular circumstances. \\
\hline
\end{tabular}
\end{table}
However, they differ in the following ways:\\


\begin{table}[H]
\centering
\caption{\textbf{Dissimilarities:}}
\renewcommand{\arraystretch}{1.3}
\begin{tabular}{|p{4cm}|p{5.5cm}|p{5.5cm}|}
\hline
\textbf{Geometric property} & \textbf{FLRW spacetime} & \textbf{LTB spacetime} \\
\hline

Symmetry & Homogeneous and isotropic cosmological model. & Spherically symmetric but radially inhomogeneous cosmological model. \\
\hline

Metric dependence & Metric depends only on cosmic time through the scale factor $a(t)$. & Metric depends on both radial coordinate $r$ and time $t$ through the function $B(t,r)$. \\
\hline

Riemann curvature tensor & Completely determined by the scale factor $a(t)$ and curvature parameter $k$. & Determined by the areal radius $B(t,r)$ and its derivatives with respect to $r$ and $t$. \\
\hline


Scalar curvature $\kappa$ & Function of cosmic time only. & Function of both radial coordinate and time. \\
\hline



Gravitational degrees of freedom & No free gravitational field; tidal forces vanish. & Tidal gravitational effects exist due to non-zero Weyl curvature. \\
\hline

Ricci soliton possibility & Ricci soliton structures exist for suitable choices of the scale factor $a(t)$. & Ricci soliton solutions may exist but depend on the radial inhomogeneity of the model. \\
\hline


Physical interpretation & Describes the large-scale homogeneous expansion of the universe. & Describes inhomogeneous cosmological evolution such as cosmic voids or gravitational collapse. \\
\hline

Curvature inheritance properties & generalized curvature Inheritance for the non-Killing vector field $\frac{\partial}{\partial t}$ and $\frac{\partial}{\partial r}.$ &  generalized curvature Inheritance for the non-Killing vector field $\frac{\partial}{\partial t}$ and $\frac{\partial}{\partial \theta}.$\\
\hline
\end{tabular}
\end{table}

\section{\bf Conclusion}
The present investigation establishes the Friedmann--Lemaître--Robertson--Walker (FLRW) spacetime as a remarkably natable geometric structure within the framework of differential geometry and relativistic cosmology. Through a systematic analysis of its curvature properties, the study demonstrates that the FLRW spacetime satisfies several important pseudosymmetric-type curvature conditions, thereby extending the understanding of its intrinsic geometric behavior beyond the classical cosmological interpretation. In particular, the relation $R \cdot R - Q(S,R)=L_CQ(g,C)$ and the associated conformal and conharmonic curvature restrictions reveal deep interactions between the curvature tensors governing the spacetime geometry.

The obtained linear dependence relations involving the Tachibana tensors \(Q(g,C)\), \(Q(S,C)\), and the tensor 
$(C \cdot R + R \cdot C)$
further highlight the highly organized curvature structure of the FLRW model. The study also proves that the spacetime belongs to the broader class of generalized Roter type manifolds and satisfies the \(\mathrm{Ein}(3)\) condition, while simultaneously exhibiting the characteristics of a 2-quasi-Einstein manifold. These findings provide a refined geometric classification of the FLRW spacetime and strengthen its significance in the study of semi-Riemannian manifolds with special curvature properties.

Although the Ricci tensor is shown to be neither cyclic parallel nor of Codazzi type, the spacetime nevertheless satisfies several notable compatibility conditions with respect to the Riemann, conformal, projective, concircular, and conharmonic curvature tensors. This demonstrates that the geometric structure of the FLRW spacetime retains a high degree of internal coherence despite the absence of stronger Ricci symmetry conditions.

The analysis of RS and CI properties provides additional insight into the dynamical and geometric evolution of the spacetime. The existence of almost Ricci soliton and \(\eta\)-Ricci Yamabe soliton structures with respect to the non-Killing vector fields $\frac{\partial}{\partial t}
\quad \text{and} \quad
\frac{\partial}{\partial r}$
indicates the presence of nontrivial geometric flows compatible with the cosmological model. In contrast, the non existence of such structures with respect to $\frac{\partial}{\partial \theta}$
emphasizes the directional dependence of these geometric phenomena. Moreover, the existence of generalized curvature inheritance for the Riemann, Weyl conformal, concircular, and conharmonic curvature tensors enriches the understanding of symmetry inheritance in cosmological spacetimes.

Finally, the comparative study between the FLRW and Lemaître--Tolman--Bondi (LTB) spacetimes illustrates both the similarities and distinctions between homogeneous-isotropic and inhomogeneous cosmological models from geometric and physical perspectives. Overall, the results obtained in this work not only deepen the mathematical understanding of the FLRW spacetime but also create new connections between differential geometry, geometric analysis, and cosmological physics, thereby opening several possible directions for future research on curvature-restricted spacetimes and their physical applications.
\section{\bf  Acknowledgment}


The Second author greatly acknowledges to The University Grants Commission, Government of India for the award of Junior Research Fellow. All the algebraic computations of Section $5-7$ are performed by a program in Wolfram Mathematica developed by the first author A. A. Shaikh.

\section{\bf Declarations}

{\bf Data availability:} Data sharing not applicable to this article as no data sets were used/generated or analyzed during the current study.\\

{\bf Competing interests:} The authors have no competing interests to declare that are relevant to the content of this article.\\

{\bf Funding:} No funding was received to assist with the preparation of this manuscript.\\


%


\end{document}